\documentclass{amsart}
\usepackage[foot]{amsaddr}
\usepackage{graphicx,epstopdf} 
\usepackage[caption=false]{subfig} 
\usepackage{tikz}
\pdfoutput=1

\newcommand{\e}{\varepsilon}
\newcommand{\edit}[1]{{\color{black}{#1}}}

\newcommand\numberthis{\addtocounter{equation}{1}\tag{\theequation}}

\newcommand{\pder}[3]{\frac{\partial^{#3} #1}{\partial #2^{#3}}}
\newcommand{\mL}{\mathcal{L}}
\newcommand{\ie}{\emph{i.e.}}

\newcommand{\cA}{\mathcal{A}}
\newcommand{\cB}{\mathcal{B}}
\newcommand{\cC}{\mathcal{C}}
\newcommand{\cD}{\mathcal{D}}

\newtheorem{thm}{Theorem}[section]

\theoremstyle{definition}
\newtheorem{definition}[thm]{Definition}

\theoremstyle{remark}
\newtheorem{remark}[thm]{Remark}

\begin{document}

\title[Absolute instabilities of travelling waves in a Keller-Segel model]{Absolute instabilities of travelling wave solutions in a Keller-Segel model}

\author{P.N. Davis$^\dagger$}
\address{$^\dagger$School of Mathematical Sciences, Queensland University of Technology, Brisbane, QLD 4000, Australia}
\author{P. van Heijster$^\dagger$}
\author{R. Marangell$^*$}
\address{$^*$School of Mathematics and Statistics, University of Sydney, Sydney, NSW 2006, Australia}

\begin{abstract}
We investigate the spectral stability of travelling wave solutions in a Keller-Segel model of bacterial chemotaxis with a logarithmic chemosensitivity function and a constant, sublinear, and linear consumption rate. Linearising around the travelling wave solutions, we locate the essential and absolute spectrum of the associated linear operators and find that all travelling wave solutions have essential spectrum in the right half plane. However, we show that in the case of constant or sublinear consumption there exists a range of parameters such that the absolute spectrum is contained in the open left half plane and the essential spectrum can thus be weighted into the open left half plane. For the constant and sublinear consumption rate models we also determine critical parameter values for which the absolute spectrum crosses into the right half plane, indicating the onset of an absolute instability of the travelling wave solution. We observe that this crossing always occurs off of the real axis.

\end{abstract}
\maketitle

\section{Introduction}
\label{S:I}
\subsection{The Keller-Segel model}\label{SS:KS}

A general Keller-Segel model of chemotaxis in one space dimension is  
\begin{equation}
\label{EQ:KSgen}
\begin{aligned}
u_t&={\e} u_{xx}-\alpha w u^m+\kappa u,\\
w_t&=\delta w_{xx}-{\beta}\left(\Phi_x(u) w\right)_x.
\end{aligned}
\end{equation}
The model represents the directed movement of a cell species $w$, such as a bacterial population, governed by the gradient of a chemical $u$. 
The function $\Phi(u)$ is the so-called chemotactic function. We take $(x,t)\in\mathbb{R}\times\mathbb{R}^+$, with $\alpha,\ \kappa\geq0, m \in \mathbb{R}$, and $\beta,\ \delta>0$ and assume that the diffusion of the chemical is taken to be much smaller than that of the bacteria, \ie\ $0\leq {\e}\ll \delta$.

Originally proposed by Keller and Segel in the 1970's (see \cite{keller1971model, keller1971traveling}) much of the focus in the literature has been on the so-called minimal Keller-Segel model (see, for example, \cite{horstmann2005boundedness,kang2007stability} and references therein, as well as the review paper \cite{horstmann20031970}). This is \eqref{EQ:KSgen} with a chemotactic function of the form $\Phi_x(u)=u$ and $\kappa=0$ (representing no growth of the chemical in the absence of the bacteria). The minimal Keller-Segel model admits solutions that blow-up in finite or infinite time \cite{horstmann20031970}.
As blow-up solutions are not biologically feasible, efforts have been made to prevent or bound blow-up solutions in the minimal Keller-Segel model by appending the model; for instance   by selecting an appropriate growth term \cite{kolokolnikov2014basic}, 
by bounding the chemotactic function \cite{horstmann2005boundedness}, \edit{or by incorporating nonlinear diffusivity \cite{wang2010chemotaxis}}.

Alternatively, by moving away from the minimal Keller-Segel model, one can find travelling wave solutions by the choice of a singular chemotactic function \cite{keller1971traveling, schwetlick2003traveling}. The literature predominantly discusses the case when the growth term $\kappa =0$, and when $\Phi(u)=\log(u)$ (see \cite{ebihara1992singular,feltham2000travelling,keller1975necessary} and the references therein). In this manuscript, we consider such a Keller-Segel model:
\begin{equation}
\label{EQ:KSalpha}
\begin{aligned}
u_{{t}}&={\e} u_{{x}{x}}-\alpha w u^m,\\
w_{{t}}&=\delta w_{{x}{x}}-{\beta}\left(\frac{wu_{{x}}}{u}\right)_{{x}}.
\end{aligned}
\end{equation}
The condition ${\beta}/\delta+m>1$ is necessary for finite solutions \cite{keller1971traveling}. It has been shown that for $m>1$ and $m<0,$ \eqref{EQ:KSalpha} admits no travelling wave solutions \cite{schwetlick2003traveling,wang2013mathematics}, thus we take $0\leq m\leq 1$. When $0\leq m \leq 1$, there are two main cases; first, for $0\leq m<1$, the model supports a travelling front of the chemical attractant coupled with a travelling pulse for the bacterial population \cite{nagai1991traveling, wang2013mathematics}. This has been used to model travelling bands of bacteria \cite{holz1978quasi, novick1984gradually}. When $m=1$, \eqref{EQ:KSalpha} supports a pair of travelling fronts and has been used to model the boundary behaviours of populations of bacteria \cite{nossal1972boundary}. See Figure \ref{fig:tw_profiles} for plots of travelling wave solutions in these two cases.

While the existence of travelling wave solutions to \eqref{EQ:KSalpha} has been studied since the model's inception, stability analysis of these travelling wave solutions has been comparatively limited.  A typical first step in the stability analysis of travelling wave solutions is to linearise around the travelling wave solution and to compute the spectrum of the resulting linearised operator. 
For travelling wave solutions in \eqref{EQ:KSalpha}, with $\e = m=0$, the essential spectrum (see Definition \ref{DEFN:EssSpect}) of the associated linear operator, dealing with instabilities at infinity, was located in \cite{nagai1991traveling}. It was shown that the essential spectrum always intersects the right half plane and so the waves are (spectrally) unstable. It is possible to shift the essential spectrum using weighted function spaces, see \S\ref{SS:Weighted_spaces}. In \cite{nagai1991traveling} a weighted function space was considered for a range of weights and it was shown that in this range the spectrum remains unstable. These results were generalised in \cite{wang2013mathematics} for $0\leq m\leq1 $. 

In this manuscript, we locate the essential spectrum associated with travelling wave solutions in \eqref{EQ:KSalpha}. 
By computing the absolute spectrum (see Definition \ref{DEFN:Abs_spect}), 
we show that for all $0\leq m<1$ there exists a range of the chemotactic parameter $\beta$, independent of the speed of the travelling wave solution, such that the essential spectrum can be weighted fully into the left half plane for an appropriate two-sided weight. See \S\ref{SS:theorem} for a more in depth explanation of the main results. 

In \S \ref{SS:SETUP}, we describe the linearised eigenvalue problem associated with a travelling wave solution to \eqref{EQ:KSalpha}, outline the particulars of spectral theory, and state our main results.  
In \S\ref{S:MZERO}, we locate the essential and absolute spectrum and explain the procedure for calculating the so-called ideal weight (see Definition \ref{DEFN:ideal_weight}), in the case of constant consumption and zero diffusivity of the attractant, \ie\ $\e=m=0$. We also calculate the range of $\beta$ values for which the essential spectrum can be weighted into the left half plane. Outside this range the travelling wave solutions are absolutely unstable.  
In \S\ref{S:SUBLIN}, we extend the results of the constant consumption case ($m=0$) to the case of sublinear $(0<m<1)$ and linear consumption ($m=1$), still in the absence of diffusion of the attractant.
While the procedures of \S\ref{S:SUBLIN} are similar to the procedures of \S\ref{S:MZERO}, the computations are algebraically more involved and therefore we split these two sections. 
In \S\ref{SS:MNONZERO_EPSNONZERO}, we include a small, non-zero, diffusivity of the attractant in the model, \ie\ $0<\e \ll1$, and show that (in)stability conditions are to leading order the same as before. 
We conclude the manuscript with a summary and discussion of future work.

\section{Set-up, definitions, and main results}
\label{SS:SETUP}
We briefly discuss the existence of travelling wave solutions to \eqref{EQ:KSalpha} and define the stability problem. Following \cite{nagai1991traveling}, we nondimensionalise \eqref{EQ:KSalpha} through the change of variables $\tilde{x}:=\sqrt{\frac{\alpha}{\delta}}{x},\ \tilde{t}:=\alpha {t}$.
Then, \eqref{EQ:KSalpha} becomes
\begin{equation}
\label{EQ:KStildex}
\begin{aligned}
u_{\tilde{t}}&=\tilde{\e} u_{\tilde{x}\tilde{x}}- w u^m,\\
w_{\tilde{t}}&=w_{\tilde{x}\tilde{x}}-\tilde{\beta}\left(\frac{wu_{\tilde{x}}}{u}\right)_{\tilde{x}},
\end{aligned}
\end{equation}
where we have set $\tilde{\e}:=\frac{{\e}}{\delta}$ and $\tilde{\beta}:=\frac{{\beta}}{\delta}$. 
We drop the tildes for notational convenience
\begin{equation}
\label{EQ:KStw1}
\begin{split}
u_t&=\e u_{xx}- w u^m,\\
w_t&=w_{xx}-\beta\left(\frac{wu_x}{u}\right)_x,
\end{split}
\end{equation}
and the conditions on our parameters are now $0\leq \e \ll1$, $\beta+m>1$ and $0\leq m\leq 1$.
\subsection{Travelling wave solutions}
\begin{figure}[t]
  \centering
\subfloat{\label{fig:tw1}\scalebox{1}{\includegraphics{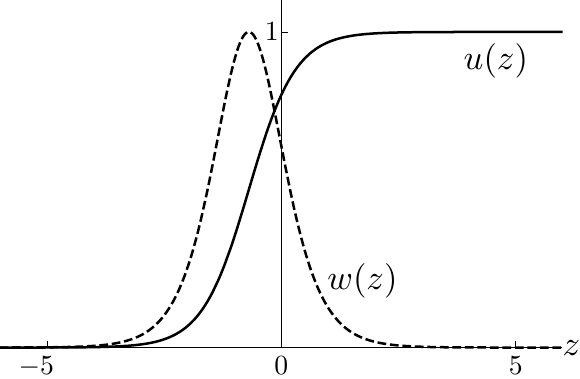}}}\hspace{.5cm}
\subfloat{\label{fig:tw2}\scalebox{1}{\includegraphics{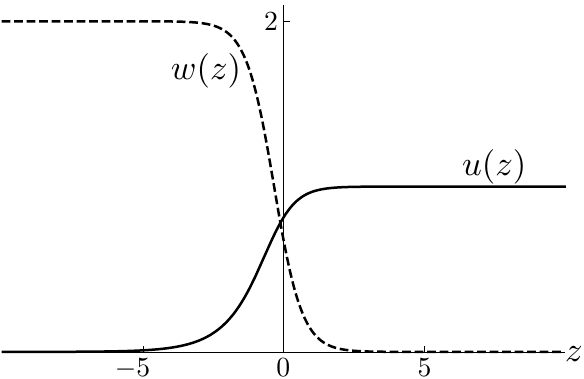}}}
  \caption{Travelling wave solutions to \eqref{EQ:KStw} for $\e=0$, $\beta=c=2$. Left panel: For $m=0$ the travelling wave solutions are a front and a pulse. Right panel: For $m=1$ the travelling wave solutions are a pair of travelling fronts.}
  \label{fig:tw_profiles}
\end{figure}
We make the change of variables $z=x-ct$, where $c>0$ is a constant, finite wave speed. In this moving frame, we have
\begin{equation}
\label{EQ:KStw}
\begin{split}
u_t&=\e u_{zz}+cu_z- w u^m,\\
w_t&=w_{zz}+cw_z-\beta\left(\frac{wu_z}{u}\right)_z.
\end{split}
\end{equation} 
Travelling wave solutions exist as stationary solutions to \eqref{EQ:KStw}, \ie\ $(u(z,t),w(z,t)) = (u(z),w(z))$ and satisfy
\begin{equation}
\label{EQ:KStw2}
\begin{split}
0 &= \e u_{zz}+cu_z- w u^m,\\
0 &= w_{zz}+cw_z-\beta\left(\frac{wu_z}{u}\right)_z\,.
\end{split}
\end{equation}
When $0\leq m<1$, travelling wave solutions satisfy \eqref{EQ:KStw2} 
with
 $$\lim_{z\to-\infty}u(z)=\lim_{z\to-\infty}w(z)=0,\quad\lim_{z\to\infty}u(z)=u_r,\quad \lim_{z\to\infty}w(z)=0,$$
where $u(z)$ is a wavefront and $w(z)$ is a pulse \cite{nagai1991traveling, wang2013mathematics} (see the left panel of Figure \ref{fig:tw_profiles}).
When $m=1$ travelling wave solutions satisfy \eqref{EQ:KStw2} with
  $$\lim_{z\to-\infty}u(z)=0,\quad \lim_{z\to-\infty}w(z)=\frac{c^2}{\beta}+\e\frac{c^2}{\beta^2} ,\quad\lim_{z\to\infty}u(z)=u_r,\quad \lim_{z\to\infty}w(z)=0,$$
 where both $u(z)$ and $w(z)$ are now wavefronts \cite{wang2013mathematics} (see the right panel of Figure \ref{fig:tw_profiles}).

Though explicit formulas for travelling wave solutions are known only for $\e = 0$ (\ie\ zero-diffusivity of the chemoattractant), the existence of travelling wave solutions in \eqref{EQ:KSalpha} has been shown for $0\leq m \leq 1$ and small enough values of the diffusivity of the chemoattractant (\ie\ $0 \leq \e \ll1$), see, for example, \cite{harley2014geometric,nagai1991traveling,wang2013mathematics} and the references therein. 
To leading order in $\e$, the profiles of travelling wave solutions are given by
\begin{align*}
u(z)&=\left( u_r^{-1/\gamma}+\sigma e^{-c(z+z^*)} \right)^{-\gamma}\,,\\
w(z)&=e^{-c(z+z^*)}\left(u(z)\right)^\beta \numberthis\label{tw_profiles_0}\,,\\
\gamma&=\frac{1}{\beta+m-1}, \quad \sigma=\frac{\beta+m-1}{c^2},
\end{align*}
where $z^*$ is a constant associated with the location of the centre of the travelling wave solution, and $u_r$ is the end state of the chemoattractant \cite{feltham2000travelling,nagai1991traveling,wang2013mathematics}. Because of translation invariance, we set $z_* = 0$, and because of scaling invariance in the nondimensionalisation of \eqref{EQ:KSalpha} to \eqref{EQ:KStildex}, we take $u_r = 1$ \cite{harley2014geometric}, in the remainder of this manuscript without loss of generality. 
Furthermore, from \cite{nagai1991traveling, wang2013mathematics} we have the following limits for the travelling wave solutions 
 \begin{align}
 \lim_{z\to-\infty} \frac{u_z}{u}=\frac{c}{\beta+m-1},\quad\quad  \lim_{z\to-\infty}\frac{w}{u^{1-m}}= \frac{c}{\beta+m-1}\left(\frac{c\e}{\beta+m-1}+c\right), \label{EQ:limits}
\end{align}
which will be useful for the stability analysis in the upcoming sections.

\subsection{The spectral problem}
To determine the stability of the travelling wave solutions $(u,w)$ of \eqref{EQ:KStw1}, we consider $U(z,t)=u(z)+p(z,t)$, and $W(z,t)=w(z)+q(z,t)$, where $p, q$ are perturbations in some appropriately chosen Banach space $\mathcal{X}$. Substituting $U$ and $W$ into \eqref{EQ:KStw} and considering only leading order terms for $p$ and $q$, we obtain the linear operator 
$\mL$ defined by,
\begin{align*}&\begin{pmatrix}
 p\\q
 \end{pmatrix}_t=\mL \begin{pmatrix}
 p\\q
 \end{pmatrix},\quad\quad
\mathcal{L}:=\begin{pmatrix}\e\partial_{zz}+c\pder{}{z}{}- m w u^{m-1}& -u^m\\ \mathcal{L}_p &\mathcal{L}_q
\end{pmatrix}\numberthis\label{EQ:Lop_gen}
\end{align*}
where
\begin{equation}
\label{EQ:Lop_gen_pq}
\begin{aligned}
\mathcal{L}_p&:=\beta\left(\frac{w_zu_z }{u^2}+\frac{wu_{zz} }{u^2}-\frac{2wu_z^2}{u^3}\right)+\beta\left(\frac{2wu_z}{u^2}-\frac{w_z}{u}\right)\pder{}{z}{}-\frac{\beta w}{u}\pder{}{z}{2},\\
\mathcal{L}_q&:=\beta\left(\frac{u_z^2}{u^2}-\frac{u_{zz}}{u}\right)+\left(c-\frac{\beta u_z}{u}\right)\pder{}{z}{}+\pder{}{z}{2}.
\end{aligned}
\end{equation}

The associated eigenvalue problem is obtained by taking perturbations of the form $\begin{pmatrix}
p(z,t)\\q(z,t)
\end{pmatrix}=e^{\lambda t}\begin{pmatrix}
p(z)\\q(z)
\end{pmatrix}$ 
where we now make the choice that $p,q \in \mathbb{H}^1(\mathbb{R})$. Here, $\mathbb{H}^1(\mathbb{R})$ is the usual Sobolev space of once (weakly) differentiable functions such that both the function and its first (weak) derivative (in $z$) are in $\mathbb{L}^2(\mathbb{R})$, \ie\ square integrable. Equation \eqref{EQ:Lop_gen} becomes 
\begin{equation}\label{EQ:eigen}
\begin{aligned}
\mL :\mathbb{H}^1(\mathbb{R})\times \mathbb{H}^1(\mathbb{R})  & \to  \mathbb{H}^1(\mathbb{R})\times \mathbb{H}^1(\mathbb{R}) \\
\mL \begin{pmatrix}
p\\q
\end{pmatrix} &= \lambda \begin{pmatrix}
p\\q
\end{pmatrix}.
\end{aligned}
\end{equation}

\subsection{Spectral stability: Background and definitions}
\label{SS:spectral_background}
A travelling wave solution is said to be spectrally stable if the spectrum of the associated linear operator $\sigma(\mL)$ is contained in the closed left half plane except for the origin. 
The spectrum $\sigma(\mL)$ is defined as follows:
\begin{definition}\label{DEFN:SPECTRUM} (\cite{sandstede2002stability} Definition 3.2) We say $\lambda\in\mathbb{C}$ is in the spectrum of a linear operator $\mathcal{L}$, denoted $\sigma(\mathcal{L})$, if the operator $\mathcal{L}-\lambda I$, where $I$ is the identity operator, is not invertible, \ie\ the inverse does not exist or is not bounded. 
\end{definition}
The spectrum of $\mL$ falls naturally into two parts, the essential spectrum, denoted $\sigma_{\rm ess}(\mathcal{L})$, and the point spectrum, denoted $\sigma_{\rm pt}(\mathcal{L})$ \cite{sandstede2000spectral}. The focus of this manuscript is on the essential spectrum of $\mL$. We refer to \S\ref{SS:POINT} for a discussion on the point spectrum of $\mL$.

\subsubsection{The essential spectrum}
We define an operator $\mathcal{T}(\lambda)$, equivalent to $\mL-\lambda I$, by transforming the eigenvalue problem into a system of first order order ordinary differential equations (ODEs);
\begin{align}
\mathcal{T}(\lambda)\textbf{p}:=\left(\frac{d}{dz}-M(z,\lambda)\right)\textbf{p}=0. \label{EQ:Top_gen_setup}
\end{align}

The essential spectrum of an operator of the form in \eqref{EQ:Top_gen_setup} is found by analysing the asymptotic behaviour of the operator $\mathcal{T}(\lambda)$. We set $M_\pm(\lambda):= \displaystyle\lim_{z\rightarrow\pm\infty}M(z,\lambda)$ and define the asymptotic operator associated with $\mathcal{T}(\lambda)$ as the piecewise constant operator
\begin{align*}
&\mathcal{T}_\infty(\lambda) := \begin{cases} \dfrac{d}{dz}-M_-(\lambda) &\mbox{if } z<0, \\[2.2mm]
\dfrac{d}{dz}-M_+(\lambda) & \mbox{if } z\ge 0. \end{cases}  \numberthis\label{EQ:Tinfty}
\end{align*}

The essential spectrum is found by analysing the dimensions of the unstable, stable and centre subspaces of $M_\pm(\lambda)$. We define the Morse index $i(A)$ of a constant matrix $A$ as the dimension of its unstable subspace, see \cite{kapitula2013spectral} Definition 3.1.9. So, for an asymptotic operator of the form of \eqref{EQ:Tinfty}, we denote the Morse indices $i_\pm:=i(M_\pm(\lambda)):=\dim(\mathbb{E}^u_\pm),$ where $\mathbb{E}^u_\pm$ denotes the unstable subspace of $M_\pm(\lambda)$ respectively.

\begin{definition} (\cite{kapitula2013spectral} Definition 3.1.11) \label{DEFN:EssSpect}   We say $\lambda \in \sigma_{\rm ess}(\mathcal{T}_\infty)$, the essential spectrum of $\mathcal{T}_\infty$, if either
   \begin{enumerate}
   \item $M_+(\lambda)\text{ and }M_-(\lambda)$ are hyperbolic with a different number of unstable matrix eigenvalues, \ie\ $i_+-i_-\neq0$; or
   \item  $M_+(\lambda) \text{ or } M_-(\lambda)$ has at least one purely imaginary matrix eigenvalue. 
   \end{enumerate}
\end{definition}

The essential spectrum is conserved under relatively compact perturbations of an operator. This follows from Weyl's essential spectrum theorem, see for example \cite{kapitula2013spectral} Theorem 2.2.6 and \cite{katoperturbation} Theorem 5.35. In a variety of operators that arise from linearisation about travelling wave solutions, including the Keller-Segel model \eqref{EQ:KStw}, the operator $\mathcal{T}_\infty$ is a relatively compact perturbation of $\mathcal{T}$ (see for example \cite{kapitula2013spectral} Theorem 3.1.11 or \cite{henry1981geometric}) and so their essential spectra coincide.

Due to the continuous dependence of $\mathcal{T}(\lambda)$ on $\lambda$ we have that the essential spectrum is bounded by the values of $\lambda$ where $M_+(\lambda) \text{ or } M_-(\lambda)$ has at least one purely imaginary matrix eigenvalue. These $\lambda$ values form curves in the complex plane referred to as the dispersion relations of the respective matrices. 

Generally, the region of the complex plane containing $\Re(\lambda)\gg1$ is not contained in the essential spectrum, \ie\ the region to the right of the essential spectrum has $i_+=i_-$. This condition is related to well-posedness of the eigenvalue problem  \cite{kapitula2013spectral} (see also the left panel of Figure~\ref{FIG:spec}) and is satisfied for the Keller-Segel model discussed in this manuscript.

\begin{remark}
 Following the terminology of \cite{kapitula2013spectral,sandstede2000spectral}, we refer to the matrix eigenvalues $\mu$ of $M_\pm(\lambda)$ as the spatial eigenvalues and to $\lambda$ as the temporal spectral parameter. Values $\lambda$ for which there is a solution to \eqref{EQ:eigen} are referred to as temporal eigenvalues. We note that temporal eigenvalues as defined here can be either in $\sigma_{\rm ess}$ or in $\sigma_{\rm pt}$. 
\end{remark}
\subsubsection{The absolute spectrum}
\label{SS:absolute}
\begin{figure}
  \centering
\scalebox{1}{\includegraphics{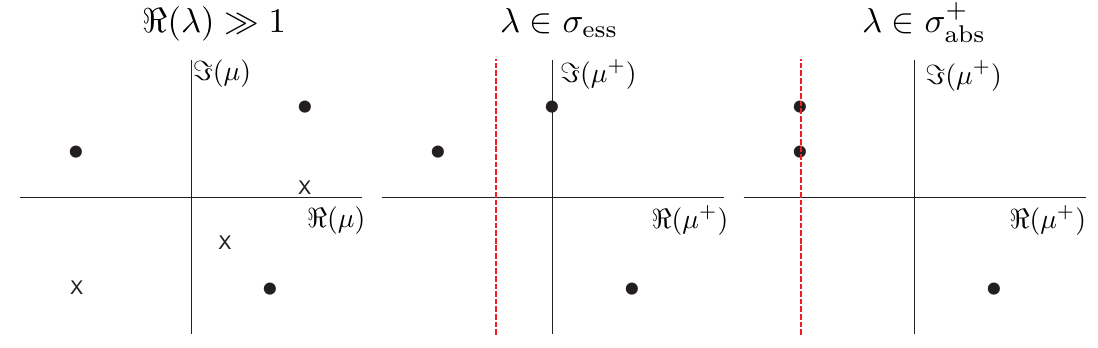}}
  \caption{A schematic of the spatial eigenvalues of the asymptotic matrices $M_+(\lambda)$ (dots) and $M_-(\lambda)$ (crosses),
with $M_\pm(\lambda)$ $3\times 3$ matrices, for three distinct values $\lambda\in\mathbb{C}$. 
Left panel: for $\Re(\lambda) \gg 1$, $M_\pm(\lambda)$ are hyperbolic and $i_\pm=2$.  
Middle panel: $\lambda\in \sigma_{\rm ess}$
since $M_+(\lambda)$ has a purely imaginary spatial eigenvalue. 
However, there exists a weight, represented by the red line, such that $i_+=2$ in this weighted space. 
So, $\lambda \notin \sigma_{\rm ess}$ in the weighted space 
(and $\lambda\notin\sigma_{\rm abs}^+$).
Right panel: $\Re(\mu_1^+)> \Re(\mu_2^+)=\Re(\mu_3^+)$, so $\lambda\in\sigma_{\rm abs}^+$ (since $i_+=2$ for $\Re(\lambda) \gg 1$, see left panel).
Observe that the order of the spatial eigenvalues persists under all weights, \ie\ the absolute spectrum does not change under weighting the space.
However, there exists a unique weight, represented by the red line, such that $\lambda$ is in the boundary of the weighted essential spectrum.
This image is adapted from Figure 3.6 of \cite{kapitula2013spectral}. 
}
  \label{FIG:spec}
\end{figure}

 The absolute spectrum, denoted $\sigma_{\rm abs}$, is not spectrum in the usual sense as it does not arise from Definition \ref{DEFN:SPECTRUM}, see, for instance, \cite{kapitula2013spectral,sandstede2002stability,sandstede2000absolute}. However, it provides important stability information as it gives an indication of how far the essential spectrum can be shifted by allowing for perturbations in weighted spaces (instead of $\mathbb{H}^1$), see also Figure~\ref{FIG:spec}. If the absolute spectrum contains values in the right half plane the solutions are said to be {\it absolutely unstable} \cite{kapitula2013spectral,sandstede2000absolute}. The absolute spectrum of $\mathcal{T}_\infty$ (equivalently of $\mathcal{T}$) is defined as follows:
\begin{definition}(\cite{sandstede2002stability} Definition 6.1)\label{DEFN:Abs_spect}
Take an $N$ dimension asymptotic operator, $\mathcal{T}_\infty$, in the form of \eqref{EQ:Tinfty}, that is well-posed in the sense that $i_+=i_-=j$ for $\Re(\lambda)\gg1$. For $\lambda\in\mathbb{C}$ we rank the $N$ spatial eigenvalues $\mu_i^\pm$ of the asymptotic matrices $M_\pm$ by the magnitude of their real parts, \ie 
\begin{align*}
\Re(\mu^\pm_1(\lambda))\geq \Re(\mu ^\pm_2(\lambda))\geq\hdots\geq \Re(\mu^\pm_j(\lambda))\geq \Re(\mu^\pm_{j+1}(\lambda))\geq\hdots\geq \Re(\mu^\pm_{N}(\lambda)).
\end{align*}
We define the sets 
\begin{align}
\sigma_{\rm abs}^+=\left\lbrace \lambda\in\mathbb{C}\left|\Re(\mu^+_j)= \Re(\mu^+_{j+1})\right.\right\rbrace\ \text{and}\ \sigma_{\rm abs}^-=\left\lbrace \lambda\in\mathbb{C}\left|\Re(\mu^-_j)= \Re(\mu^-_{j+1})\right.\right\rbrace,
\end{align}
and the absolute spectrum of $\mathcal{T}_\infty$ (and of $\mathcal{T}$) is $\sigma_{\rm abs} :=\sigma_{\rm abs}^+\cup\sigma_{\rm abs}^-$. 
 \end{definition}
Due to the continuous dependence of $\mathcal{T}$ on $\lambda$, the Morse indices will only change upon crossing one of the dispersion relations and so the absolute spectrum will always be to the left of the rightmost boundary of the essential spectrum. That is, moving $\lambda$ from right to left in the complex plane we will first encounter a dispersion relation of either $M_\pm(\lambda)$ before (potentially) encountering absolute spectrum, see also Figure~\ref{FIG:spec}.

\begin{remark}
\label{R:GENABS}
For an operator $\mathcal{T}$, with Morse indices $i_+=i_-=j$ in the region to the right of the essential spectrum, the set of $\lambda\in\mathbb{C}$ with $\Re(\mu^+_i(\lambda))= \Re(\mu^+_{i+1}(\lambda))$ or $\Re(\mu^-_i(\lambda))= \Re(\mu^-_{i+1}(\lambda))$ where $i\neq j$ is referred to as the {\it generalised absolute spectrum}. 
\end{remark}
\subsubsection{Weighted spaces}
\label{SS:Weighted_spaces}

The presence of essential spectrum of a linear operator in the right half plane implies instability of the travelling wave solution in $\mathbb{H}^1$. However, for many travelling wave solutions that are widely considered `stable', the linearised operator associated with them has essential spectrum in the right half plane; one such example is the well-known Fisher-Kolmogorov-Petrovsky-Piscounov (F-KPP) equation. A resolution proposed for this apparent contradiction is to work in an appropriately weighted space \cite{sattinger1977weighted}. Weighting the space adjusts the types of perturbations allowed. 
Following \cite{kapitula2013spectral}, we define
the weighted space $\mathbb{H}^1_\nu(\mathbb{R})$ by the norm
\begin{align}
\|p\|_{\mathbb{H}^1_\nu} =\|e^{\nu z}p\|_{\mathbb{H}^1}=\|\tilde{p}\|_{\mathbb{H}^1},\label{EQ:weighted_norm}
\end{align}
where $\tilde{p}:=e^{\nu z} p$. So, $p\in \mathbb{H}^1_\nu$ if and only if $\tilde{p}\in \mathbb{H}^1$. We define $\mathbb{L}^2_\nu$ similarly. 
The weight provides information as to whether the travelling wave solutions are more sensitive to perturbations in front of the wavefront (\ie\ as $z\to\infty$) or behind the wavefront (\ie\ as $z\to-\infty$). In other words, if $\nu>0$ then the perturbation $p(z,t)$ must decay at a rate faster than $e^{-\nu z}$ as $z\to\infty$, while it is allowed to grow exponentially at any rate less than $e^{-\nu z}$ as $z\to-\infty$. 
We can also consider a two-sided weight
\begin{equation}\label{EQ:two_sided_weight}\nu=\begin{cases} \nu_- \mbox{ if} & z \leq0, \\ \nu_+\mbox{ if} & z>0,\end{cases} \end{equation}
which forces the perturbation to decay exponentially in both directions. It turns out that we need to consider a two-sided weight \eqref{EQ:two_sided_weight} in the case of the Keller-Segel model \eqref{EQ:KStw1}.

A practical consequence of considering $\mathcal{L}$ on weighted function spaces is that the essential spectrum is moved. 
In particular, 
assume we have an operator $\mathcal{T}$ of the form of \eqref{EQ:Top_gen_setup} coming from the linearisation around a travelling wave solution and with asymptotic operator \eqref{EQ:Tinfty}. 
The operator $\mathcal{T}(\lambda)$ in the weighted space is given by
$$\mathcal{T}(\lambda)\tilde{\textbf{p}}=\textbf{p}'-(M(z,\lambda)+\nu I)\textbf{p}=0,$$
with asymptotic matrices $M_\pm(\lambda)+\nu I$ \cite{kapitula2013spectral}.
So, we need to consider the magnitude and sign of the real part of the spatial eigenvalues compared to the weight, \ie\ we consider $\mu-\nu$, the spatial eigenvalues of $M_\pm(\lambda)+\nu I$, instead of $\mu$, the spatial eigenvalues of $M_\pm(\lambda)$. See Figure~\ref{FIG:spec}.
If the operator $\mathcal{T}$ has essential spectrum in the right half plane in the unweighted space, 
weights of interest are those that move this essential spectrum into the open left half plane.
If such weights $\nu$ exist (and if there is no point spectrum in the right half plane), we say the travelling wave solution is spectrally stable in $\mathbb{H}^1_\nu(\mathbb{R})$ and
it is referred to as being {\it transiently unstable} \cite{sandstede2000absolute, sherratt2014mathematical}.

Since the order of the spatial eigenvalues is not changed,
the absolute spectrum is unaffected by weighting the function space and
the presence of absolute spectrum in the right half plane indicates an {\it absolute instability}. 
In particular, in the case of an absolute instability no weights can be found that move the essential spectrum into the left half plane since the absolute spectrum is to the left of the rightmost boundary of the essential spectrum.

\subsection{Main results}
\label{SS:theorem}
In this section, we state the main results of this manuscript related to the location of the absolute spectrum of travelling wave solutions supported by \eqref{EQ:KStw1}.
\begin{thm}
\label{TH:MAIN1}
Assume that $c>0, 0 \leq m < 1$ and $\beta> 1-m$.
Let $\beta_{\rm crit}$ be the unique real root larger than one of 
\begin{equation}
\begin{aligned}
f(\beta) =\ &310 \beta^{10} - 3234 \beta^9 + 17112 \beta^8 - 49101 \beta^7 + 76180 \beta^6 - 58398 \beta^5 \\
& + 10056 \beta^4+ 15040 \beta^3 - 9680 \beta^2 + 1716 \beta -4. \label{EQ:10th_order_poly}
\end{aligned}
\end{equation}
Then, there exists an $\e_0>0$ such that for all $0\leq \e <\e_0$ the absolute spectrum of $\mathcal{L}$ given in \eqref{EQ:Lop_gen} is fully contained in the left half plane for all $1-m<\beta<\beta_{\rm crit}^m(\e)$, with $\beta_{\rm crit}^m(\e)$ to leading order given by $\beta_{\rm crit}^m := \beta_{\rm crit}(1-m)$.
Crucially, at $\beta=\beta_{\rm crit}^m(\e)$ the absolute spectrum crosses into the right half plane off of the real axis with increasing $\beta$. For $\beta>\beta_{\rm crit}^m(\e)$ the absolute spectrum of $\mathcal{L}$ \eqref{EQ:Lop_gen} contains values in the right half plane and the travelling wave solutions of \eqref{EQ:KStw1} are thus absolutely unstable.

For $m=1$, the absolute spectrum of $\mathcal{L}$ \eqref{EQ:Lop_gen} includes the origin for all parameter values. 
\end{thm}

The fact that the polynomial $f$ \eqref{EQ:10th_order_poly} has only one real root larger than one follows directly from Sturm's Theorem, see, for instance, Theorem 6.3d in \cite{henrici1978applied}. In particular, $\beta_{\rm crit} \approx 1.6195$.
Moreover, for every $0\leq m<1$ and $1<\beta<\beta_{{\rm  crit}}^m(\e)$ there exists a range of two-sided weights $\nu$ \eqref{EQ:two_sided_weight} such that weighted essential spectrum is contained in the open left half plane, see Remark~\ref{REM:weighted_dispersion_range} and Remark~\ref{REM:weighted_dispersion_range2}. Also, observe that the above leading order results are independent of the wave speed $c$, see Remark~\ref{REM:scaling_invariance_c}.

So, we fully classify the (in)stabilities coming from the weighted essential spectrum of travelling wave solutions of \eqref{EQ:KStw1} for the complete parameter range for which travelling wave solutions exist, \ie\ for $0\leq m \leq 1$ and $1-m<\beta$ \cite{schwetlick2003traveling,wang2013mathematics}. In essence, we obtain the complete picture of the essential spectrum, extending the initial results obtained in \cite{nagai1991traveling,wang2013mathematics}.

As we are primarily concerned with the absolute spectrum, we define the {\it ideal weight} as the weight such that the weighted dispersion relations intersect the rightmost points of the absolute spectrum.
\begin{definition}\label{DEFN:ideal_weight}
The ideal weight for the operator \eqref{EQ:Lop_gen} is the unique two-sided weight such that the dispersion relations of $M_\pm (\lambda)+\nu_\pm I$ intersect the leading edges of the $\sigma_{\rm abs}^\pm$ respectively.
\end{definition}
 This definition is motivated by the fact that as $\beta$ increases, the ideally weighted essential spectrum and the absolute spectrum cross into the right half plane simultaneously. 
 
\section{Constant consumption and zero diffusivity of the chemoattractant}
\label{S:MZERO}
For clarity of presentation, we first prove Theorem \ref{TH:MAIN1} in the case of constant consumption ($m=0$) and zero diffusivity of the chemoattractant ($\e=0$). We show that the absolute spectrum is contained in the left half plane when $1<\beta<\beta_{\rm crit}$ (with $\beta_{\rm crit}$ the root of \eqref{EQ:10th_order_poly}), while it contains values in the right half plane when $\beta>\beta_{\rm crit}$. Consequently, when $1<\beta<\beta_{\rm crit}$, there exists a two-sided weight $\nu$ \eqref{EQ:two_sided_weight} such that the essential spectrum is contained in the open left half plane in the ideally weighted space,
while all travelling wave solutions are absolutely unstable when $\beta\geq\beta_{\rm crit}$.
\subsection{Set-up}
\label{SS:SETUP2}
In the $\e=m=0$ case, the eigenvalue problem \eqref{EQ:eigen} reduces to
 \begin{align*}
\mL\begin{pmatrix}p\\ q
\end{pmatrix}=\lambda\begin{pmatrix}p\\ q
\end{pmatrix},\quad\text{with}\quad \mathcal{L}&=\begin{pmatrix}c\pder{}{z}{}& -1\\ \mathcal{L}_p &\mathcal{L}_q
\end{pmatrix},
\end{align*}
where $\mL_p$ and $\mL_q$ are given by \eqref{EQ:Lop_gen_pq}, restated here for convenience,
\begin{subequations}
\begin{align*}
\mathcal{L}_p&:=\beta\left(\frac{w_z u_z }{u^2}+\frac{wu_{zz} }{u^2}-\frac{2wu_z^2}{u^3}\right)+\beta\left(\frac{2wu_z}{u^2}-\frac{w_z}{u}\right)\pder{}{z}{}-\frac{\beta w}{u}\pder{}{z}{2},\\
\mathcal{L}_q&:=\beta\left(\frac{u_z^2}{u^2}-\frac{u_{zz}}{u}\right)+\left(c-\frac{\beta u_z}{u}\right)\pder{}{z}{}+\pder{}{z}{2}.
\end{align*}
\end{subequations}
Here $(u,w)$ are the (explicit) travelling wave solutions given in \eqref{tw_profiles_0}. We define the operator $\mathcal{T}_0(\lambda)$, equivalent to $\mL-\lambda I$, by setting $s=q_z$. The operator $\mathcal{T}_0(\lambda)$, with $p,q\in \mathbb{H}^1(\mathbb{R})$ and $s\in \mathbb{L}^2(\mathbb{R})$, is given by
\begin{align*}
&\mathcal{T}_0(\lambda)\begin{pmatrix}
p\\q\\s
\end{pmatrix}:=\begin{pmatrix}
p\\q\\s
\end{pmatrix}'-M_0(z,\lambda)\begin{pmatrix}
p\\q\\s
\end{pmatrix}=0,
&&M_0(z,\lambda):=\begin{pmatrix}
 \frac{\lambda }{c} & \frac{1}{c} & 0 \\
 0 & 0 & 1 \\
 \mathcal{A}_0 &  \mathcal{B}_0 &  \mathcal{C}_0
\end{pmatrix},\numberthis\label{EQ:T0}
\end{align*}
with
\begin{align*}
\mathcal{A}_0&=\beta\left(\frac{2wu_{z}^2}{u^3}-\frac{w_zu_{z} }{u^2}-\frac{wu_{zz} }{u^2}\right)+\frac{\lambda \beta}{c}\left(\frac{w_z}{u}-\frac{2wu_{z}}{u^2}\right)+\frac{\lambda^2\beta w}{c^2 u},\\
\mathcal{B}_0&=\beta\left(\frac{u_{zz}}{u}-\frac{u_{z}^2}{u^2}\right)+\frac{\beta}{c}\left(\frac{w_z}{u}-\frac{2wu_{z}}{u^2}\right)+\frac{\lambda\beta}{c^2}\left(\frac{w}{u}\right)+\lambda,\\
\mathcal{C}_0&=\frac{\beta u_{z}}{u}-c+\frac{\beta}{c}\frac{w}{u}.
\end{align*}
\subsection{Essential spectrum}
\label{SS:ESS}
We first locate the essential spectrum in the unweighted function space. We calculate the dispersion relations of the asymptotic matrices as these act as the boundaries of the essential spectrum. 
From \eqref{EQ:KStw2}, with $\e=0$, we have $u_z=w/c$ and by integrating the second equation we get $w_z=-cw+\beta\left(wu_z/u\right)$ (where the integration constant is zero \cite{harley2014geometric,keller1971traveling}).
Thus, all terms of $M_0$ can be written in terms of $w/u$.
From \eqref{EQ:limits}, or directly from the travelling wave profiles \eqref{tw_profiles_0}, we have, 
\begin{align*}
  \lim_{z\to\infty} \frac{w}{u}=0,\quad 
  \lim_{z\to-\infty} \frac{w}{u}=\frac{c^2}{\beta-1}
\end{align*}
Using these facts, the limits of $\mathcal{A}_0$, $\mathcal{B}_0$ and $\mathcal{C}_0$ as $z\rightarrow\pm\infty$, denoted $\mathcal{A}_{0}^{\pm}$, $\mathcal{B}_{0}^{\pm}$ and $\mathcal{C}_{0}^{\pm}$, are straightforward to compute and are, respectively, given by
\begin{align*}
\mathcal{A}_{0}^+&=0, & 
 \mathcal{A}_{0}^-&=\frac{\beta  \lambda  \left((\beta -1) \lambda -c^2\right)}{(\beta -1)^2}, \\
\mathcal{B}_{0}^+&=\lambda, & 
 \mathcal{B}_{0}^-&= \frac{\left(2\beta ^2-3 \beta +1\right) \lambda -c^2 \beta }{(\beta -1)^2}, \\
 \mathcal{C}_{0}^+&=  -c,&
 \mathcal{C}_{0}^-&=  \frac{c (\beta +1)}{\beta -1}.
\end{align*}
We also define the asymptotic matrices,
\begin{align}
 M_0^\pm(\lambda):=\displaystyle\lim_{z\rightarrow\pm\infty}M_0(z,\lambda)=\begin{pmatrix}
 \frac{\lambda }{c} & \frac{1}{c} & 0 \\
 0 & 0 & 1 \\
 \mathcal{A}_0^\pm &  \mathcal{B}_0^\pm &  \mathcal{C}_0^\pm
\end{pmatrix}.
\end{align}
The dispersion relations of $M_0^+$ are
\begin{align}
&\lambda=-k^2+ick,\quad\text{and}\quad\lambda=ick, \label{EQ:dispersion_plus}
\end{align}
where $k\in\mathbb{R}$ and where $\mu=ik$ is a purely imaginary spatial eigenvalue of $M_0^+$. Note that the imaginary axis is one of the dispersion relations, while the other is a parabola opening to the left half plan with vertex at the origin.

The dispersion relations of $M_0^-$ are given by
\begin{align}
\lambda ^2+\left(k^2-\frac{i (\beta -2) c k}{\beta -1}\right)\lambda+\frac{(\beta +1) c^2 k^2}{\beta -1}+i c k\left(\frac{ \beta  c^2}{(\beta -1)^2}- k^2\right)=0,  \label{EQ:dispersion_minus}
\end{align}
where $k\in\mathbb{R}$ and where $\mu=ik$ is a purely imaginary spatial eigenvalue of $M_0^-$. Equation \eqref{EQ:dispersion_minus} is quadratic in the temporal parameter $\lambda$ and cubic in the parameter $k$ (and thus in the spatial eigenvalue).
\begin{figure}
  \centering
\scalebox{1}{  \includegraphics{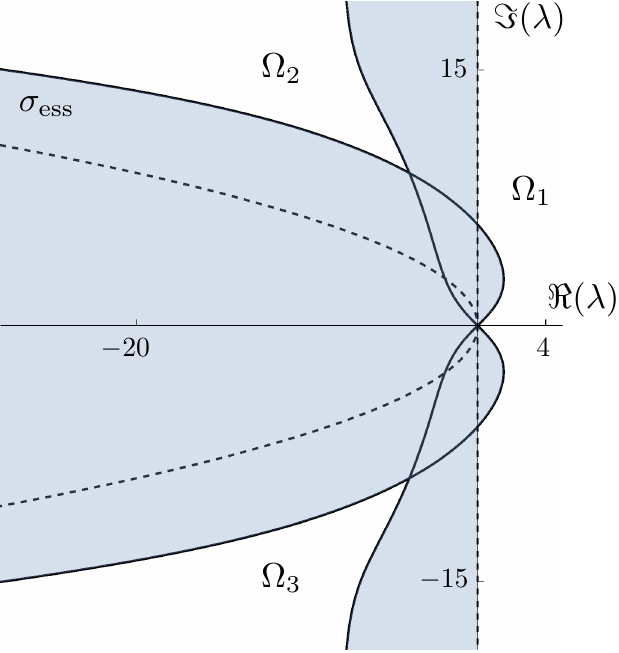}}
  \caption{The essential spectrum $\sigma_{\rm ess}$ of the operator $\mL$ about the travelling wave solutions $\left(u, w\right)$ \eqref{tw_profiles_0} for $\e=m=0$ and $\beta=c=2$. The solid curves are the dispersion relations of $M_0^-$, while the dashed curves are the dispersion relations of $M_0^+$. The shaded region is $\lambda\in\mathbb{C}$ such that $i_+\neq i_-$ and the essential spectrum is the union of the shaded region and the dispersion relations. Observe that the entire imaginary axis is included in the essential spectrum. The general shape of the unweighted essential spectrum is qualitatively similar for all values $\beta>1$, while changing the wave speed $c$ only affects the scaling of the image, see Remark \ref{REM:scaling_invariance_c}. Note this figure is a slight correction to Figure 6 from \cite{harley2015numerical}.}\label{FIG:ess_spect}
\end{figure}

The boundary of the essential spectrum is traced out by the solutions $\lambda\in\mathbb{C}$, parametrised by $k$, from \eqref{EQ:dispersion_plus} and \eqref{EQ:dispersion_minus}. We label the connected set containing $\Re(\lambda)\gg1$ as $\Omega_1$, see Figure \ref{FIG:ess_spect}. For $\lambda\in\Omega_1$, we have that the dimensions of the unstable subspaces of $M_0^\pm$ are both two, \ie\ $i_\pm=2$. There are two other regions in the complex plane where $i_+=i_-$. We denote these regions $\Omega_2$ and $\Omega_3$, see Figure \ref{FIG:ess_spect}. The remaining part of the complex plane is the essential spectrum. It is clear from Figure \ref{FIG:ess_spect} that part of the essential spectrum is in the right half plane. This agrees with previous results; by considering \eqref{EQ:dispersion_minus} for small $|k|$ values it was shown all travelling wave solutions for $\e=m=0$ are unstable in the unweighted space \cite{nagai1991traveling}.

\subsection{The weighted essential spectrum and the absolute spectrum}
\label{S:Weights_eps_zero_m_zero}
To further investigate the stability properties of the travelling wave solutions, we consider the spectrum in various two-sided weighted spaces, locate the absolute spectrum and identify the ideal weight. We substitute $\tilde{\textbf{p}}=e^{\nu z} \textbf{p}$, where $\textbf{p}=(p,q,s)^T$, into \eqref{EQ:T0} and consider the weighted space $\mathbb{H}^1_\nu$ \eqref{EQ:weighted_norm} with $\nu$ a two-sided weight \eqref{EQ:two_sided_weight}. This substitution transforms \eqref{EQ:T0} into
$$\mathcal{T}_0(\lambda)\tilde{\textbf{p}}=\textbf{p}'-(M_0(z,\lambda)+\nu I)\textbf{p}=0,$$
with $M_0(z,\lambda)$ as given in \eqref{EQ:T0}.
The essential spectrum in the weighted space is bounded by the dispersion relations of the asymptotic matrices $M_0^\pm+\nu_\pm I$.

\subsubsection{The weighted dispersion relations and absolute spectrum from $M_0^+$}

First, we consider the dispersion relations of $M_0^+(\lambda)+\nu_+ I$;
\begin{align*}
\lambda=-c \nu_+ +i c k,\quad \text{and}\quad \lambda=-k^2-\nu_+(c-\nu_+)+i (c k-2 k \nu_+ ).\numberthis\label{EQ:dispersion_plus_weighted}
\end{align*}
For $\nu_+\in (0,c)$ the real part of the dispersion relations \eqref{EQ:dispersion_plus_weighted} have strictly negative real parts and the furthest left these relations can be shifted is for the ideal weight $\nu_+^*=c/2$. Under this weight, the dispersion relations \eqref{EQ:dispersion_plus_weighted} reduce to
\begin{align}
\lambda=-\frac{c^2}{2} +i c k,\quad \text{and}\quad \lambda=-\frac{c^2}{4}-k^2.\label{EQ:dispersion_plus_max_weighted}
\end{align}

\begin{figure}
  \centering
\scalebox{.9}{\subfloat{\label{fig:a}\includegraphics[trim={.6cm 0 0 0}, clip]{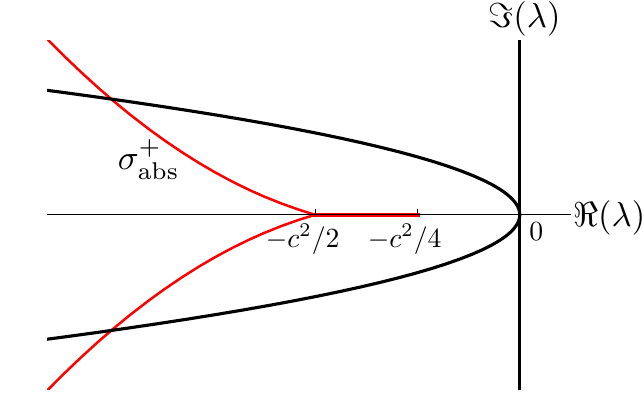}}}\hspace{.1cm}
\scalebox{.9}{\subfloat{\label{fig:b}\includegraphics[trim={.4cm 0 0 0}, clip]{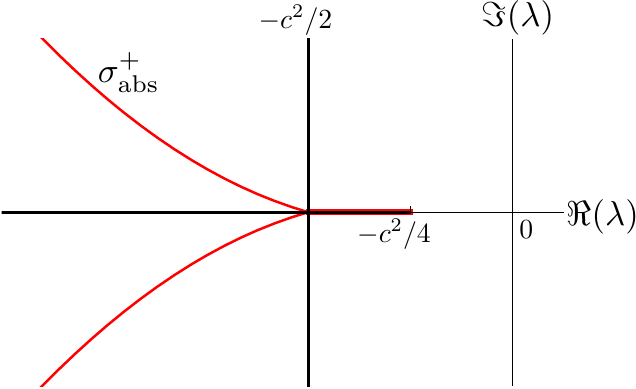}}}
  \caption{The subset of the absolute spectrum $\sigma_{\rm abs}^+$ (red) and the dispersion relations of $M_0^++\nu_+ I$ (black). Left panel: the dispersion relations \eqref{EQ:dispersion_plus_weighted} in the unweighted space, \ie\ $\nu_+=0$. The imaginary axis is one of the dispersion relations. Right panel: the ideally weighted dispersion relations \eqref{EQ:dispersion_plus_max_weighted}, \ie\ $\nu_+^*=c/2$. Note that the parabola from the left panel collapses to the real line under the ideal weight.}
  \label{fig:absplus}
\end{figure}

Next, we calculate $\sigma_{\rm abs}^+$, the subset of the absolute spectrum arising from the spatial eigenvalues for $z\to\infty$.
Since $i_+=2=i_-$ for $\Re(\lambda) \gg1$, we search for $\lambda\in\mathbb{C}$ such that the spatial eigenvalues with the second and third largest real part have the same real part (see Definition \ref{DEFN:Abs_spect}).  
The spatial eigenvalues of $M_0^+$ are
\begin{align}
\mu_{1}^+=\frac{\lambda}{c},\quad\mu_{2}^+=\frac{-c+\sqrt{c^2+4\lambda}}{2}, \quad\mu_3^+=\frac{-c-\sqrt{c^2+4\lambda}}{2}.\label{EQ:spatial_eigenvalues_plus}
\end{align}
 For $\Re(\lambda)\geq-\frac{c^2}{2}$, we have that $\Re(\mu^+_1)\geq\Re(\mu^+_2)\geq \Re(\mu^+_3)$. So, the absolute spectrum in this region is given by $\lambda\in\mathbb{C}$ such that $\Re(\mu^+_2)=\Re(\mu^+_3)$. That is, $\left\lbrace\lambda\in\mathbb{R}\left|-\frac{c^2}{2}\leq\lambda\leq\frac{-c^2}{4}\right.\right\rbrace$.
For $\Re(\lambda)<-\frac{c^2}{2}$, we have that $\mu_2^+$ has the largest real part and the absolute spectrum in this region is thus given by $\lambda\in\mathbb{C}$ such that $\Re(\mu_1)=\Re(\mu_3)$. That is, $\left\lbrace \lambda=\lambda_1+i\lambda_2,\ \lambda_1,\lambda_2\in\mathbb{R}\left| \lambda_1<-\frac{c^2}{2};\ \lambda_2=\pm\lambda_1\left(1+\frac{2 \lambda_1}{c^2}\right)\right. \right\rbrace$.
So, $\sigma_{\rm abs}^+$ is given by
\begin{equation}\label{EQ:abs_plus_m0}
   \begin{aligned}
\sigma_{\rm abs}^{+}&=\left\lbrace\lambda\in\mathbb{R}\left|-\frac{c^2}{2}\leq\lambda\leq\frac{-c^2}{4}\right.\right\rbrace\cup\\
&\left\lbrace \lambda=\lambda_1+i\lambda_2,\ \lambda_1,\lambda_2\in\mathbb{R}\left| \lambda_1<-\frac{c^2}{2};\ \lambda_2=\pm\lambda_1\left(1+\frac{2 \lambda_1}{c^2}\right)\right. \right\rbrace.
\end{aligned}
\end{equation}
Obviously, $\sigma_{\rm abs}^+$ is fully contained in the left half plane for all $c>0$. Consequently, no absolute instabilities arise from $z\to\infty$. See Figure \ref{fig:absplus} for a plot of $\sigma_{\rm abs}^+$ \eqref{EQ:abs_plus_m0} and the ideally weighted dispersion relations \eqref{EQ:dispersion_plus_max_weighted} and the unweighted dispersion relations \eqref{EQ:dispersion_plus} (or  \eqref{EQ:dispersion_plus_weighted} with $\nu_+=0$).

\subsubsection{The weighted dispersion relations and absolute spectrum from $M_0^-$}
The characteristic equation of $M_0^-$ is given by 
\begin{align*}
\mu^3-\mu^2\left(\frac{(\beta +1) c}{\beta -1}+\frac{\lambda }{c}\right)+\mu\left(\frac{(2-\beta ) \lambda }{\beta -1}+\frac{\beta  c^2}{(\beta -1)^2}\right)+\frac{\lambda ^2}{c}=0,\numberthis \label{EQ:char_poly_minus}
   \end{align*} 
and the dispersion relations of $M_0^-+\nu_- I$ are implicitly given by 
  \begin{equation}\label{EQ:dispersion_weighted_minus}
 \begin{aligned}
 &\lambda^2+\left(\frac{c(2-\beta ) (i k-\nu_-)}{\beta -1}-(i k-\nu_-)^2\right)\lambda+\frac{\beta  c^3 (i k-\nu_-)}{(\beta -1)^2}\\
 &\qquad-\frac{(\beta +1) c^2 (i k-\nu_-)^2}{\beta -1}+ c(i k-\nu_-)^3=0.
 \end{aligned}
 \end{equation}

\begin{figure}
 \begin{center}
\scalebox{1}{ \includegraphics{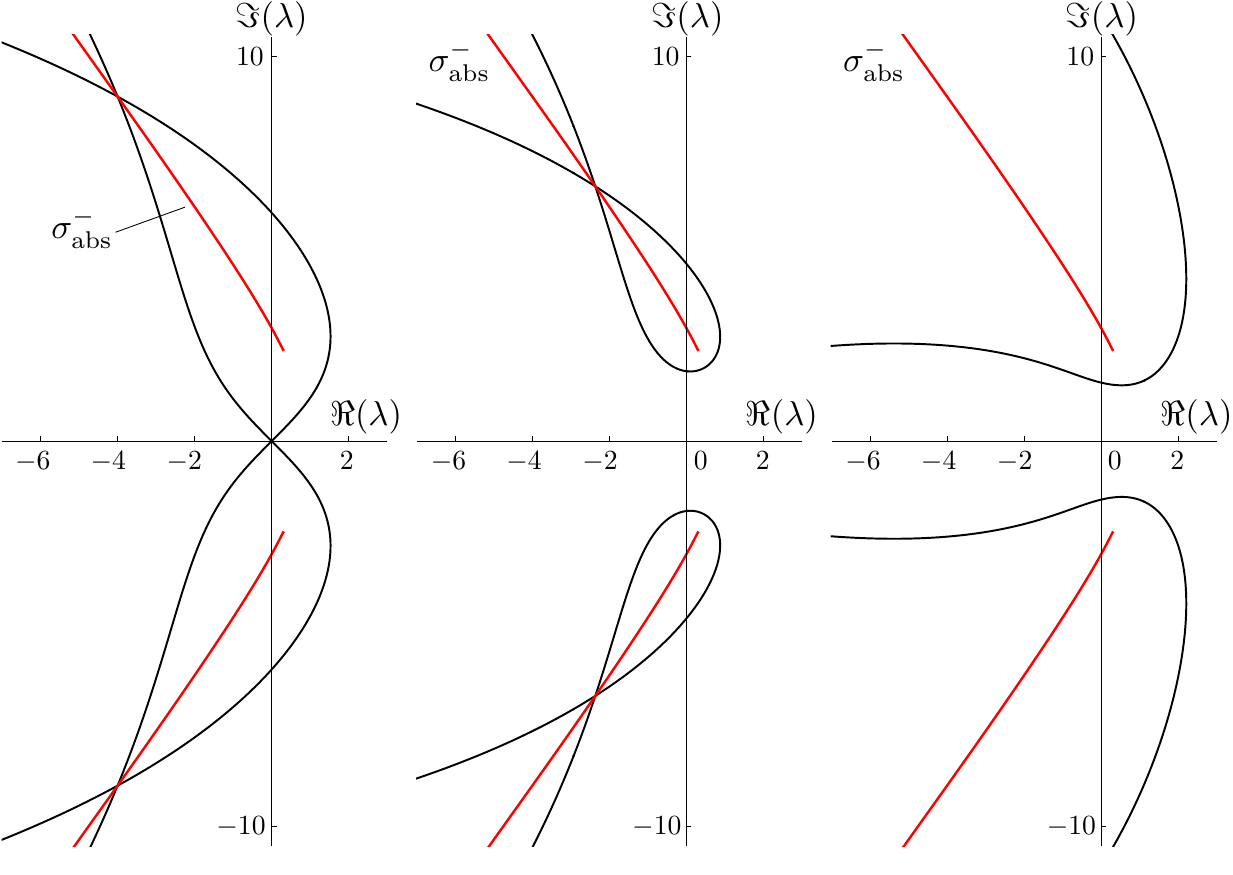}}
 \caption{
 The subset of the absolute spectrum $\sigma_{\rm abs}^-$ (red) and the dispersion relations of $M_0^-+\nu_- I$ (black) for $\beta=c=2$ and various weights $\nu_-$.  
The dispersion relations \eqref{EQ:dispersion_weighted_minus} in the unweighted space (left panel), a weighted space with $\nu_-=-1/4$ (middle panel), and a weighted space with $\nu_-=-3/2$ (right panel).
As $\nu_-$ is further decreased, the dispersion relations move further into the right half plane. For $\nu_->0$, the leading edge of the weighted dispersion relation also moves further into the right half plane.  }
 \label{FIG:various_weighted_minus}
 \end{center}
 \end{figure}

For a fixed $\beta$ and $c$ and for various weights $\nu_-$, we can plot the weighted dispersion relations \eqref{EQ:dispersion_weighted_minus}, see, for example, Figure \ref{FIG:various_weighted_minus}. 
Observe that the weighted dispersion relations \eqref{EQ:dispersion_weighted_minus} have self-intersections for some $\lambda\in\mathbb{C}$ over a large range of weights $\nu_-$, including $\nu_-=0$ (related to the unweighted space). This self-intersection corresponds to two complex roots of the characteristic polynomial \eqref{EQ:char_poly_minus} of the form $\mu_{1,2}^-=-\nu_-+ik_{1,2}$ with $k_{1,2}\in\mathbb{R}$. Thus, we have $\Re(\mu_1^-)=\Re(\mu_2^-)$, while the third spatial eigenvalue $\mu_3^-$ has a larger real part. Consequently, the $\lambda$ value at the self-intersection is part of the absolute spectrum. 

\begin{figure}
 \begin{center}
 \scalebox{.9}{\includegraphics[trim={0 .5cm 0 0}, clip]{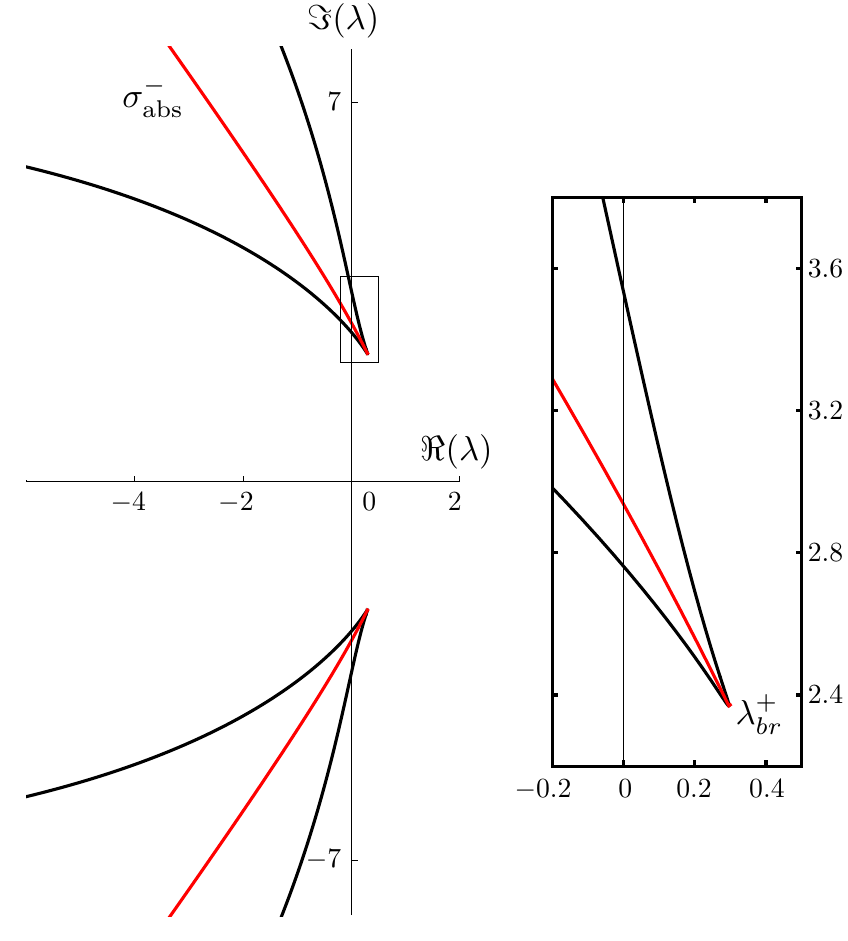}}
 \caption{
 The subset of the absolute spectrum $\sigma_{\rm abs}^-$ (red) and the ideally weighted dispersion relations of $M_0^-+\nu_-^* I$ (black) for $\beta=c=2$, where the ideal weight  $\nu_-^*\approx-0.73$.
The weighted dispersion relations form cusps whose tips coincide with the leading edge of the absolute spectrum, \ie\ the branch points $\lambda_{br}^\pm$ \edit{(see Remark \ref{REM:branch_points})}. Since the absolute spectrum, and thus the essential spectrum, enter into the right half plane, the travelling wave solution is absolutely unstable for this parameter set.}
 \label{FIG:critical_weight}
 \end{center}
 \end{figure}

There exists some weight $\nu_-^*<0$ such that the self-intersection vanishes for $\nu_-<\nu_-^*$, see, for instance, the right panel of Figure~\ref{FIG:various_weighted_minus}. 
For $\nu_-=\nu_-^*$, the self-intersection forms a cusp of the weighted dispersion relations \eqref{EQ:dispersion_weighted_minus} and is thus 
the ideal weight, see Figure~\ref{FIG:critical_weight}.
For $\nu_->\nu_-^*$, the self-intersections trace out the subset of the absolute spectrum $\sigma_{\rm abs}^-$. This allows us to directly locate $\sigma_{\rm abs}^-$ using a find root procedure on the dispersion relations of $M_0^-+\nu_-I$. 
Values $\lambda\in\sigma_{\rm abs}^-$ such that there is a second order root (in $\mu$) of the characteristic polynomial \eqref{EQ:char_poly_minus} are referred to as {\it branch points} $\lambda_{br}$, see Remark \ref{REM:branch_points} and Figure~\ref{FIG:critical_weight}. 
For the Keller-Segel model, the cusp of the ideally weighted dispersion relations corresponds to the second order root  and so the branch points are the rightmost points of $\sigma_{\rm abs}^-$, see Figure~\ref{FIG:critical_weight}. 

To locate the branch points $\lambda_{br}$, we treat the characteristic polynomial \eqref{EQ:char_poly_minus} as a cubic polynomial in $\mu$ and determine the second order roots. This boils down to finding $\lambda\in\mathbb{C}$ such that the discriminant of \eqref{EQ:char_poly_minus} is zero. That is, we solve
\begin{align*}
&\lambda^5+\frac{(2
   \beta -1)^2 c^2 \lambda ^4}{4 (\beta -1)^2}+\frac{\beta(18 \beta^2 -37\beta+20) c^4 \lambda ^3}{2 (\beta -1)^3}+\frac{\beta(5 \beta^3 -28\beta^2+50\beta-26) c^6 \lambda ^2}{4 (\beta
   -1)^4}\\
   &\quad\quad-\frac{\beta(\beta^2 -6\beta +2  ) c^8 \lambda }{2
   (\beta -1)^4}+\frac{\beta ^2 c^{10}}{4 (\beta -1)^4}=0. \numberthis \label{EQ:discriminant}
   \end{align*}
 We look for roots of \eqref{EQ:discriminant} that correspond to the two smallest spatial eigenvalues having the same real part, \ie\ the values $\lambda\in\sigma_{\rm abs}^-$ that solve \eqref{EQ:discriminant}.
For given parameters, we find a pair of complex conjugate solutions to \eqref{EQ:discriminant} that are in the absolute spectrum; these solutions are the branch points $\lambda_{br}^\pm$ that form the leading edge of $\sigma_{\rm abs}^-$. Note that the other three roots of \eqref{EQ:discriminant} are part of the generalised absolute spectrum, see Remark~\ref{R:GENABS}.

Locating the branch points $\lambda_{br}^\pm$ also allows us to compute the ideal weight $\nu_-^*$, since $\nu_-^*$ corresponds to the negative of the real part of the second order root $\mu$ evaluated at the branch point $\lambda_{br}^\pm$. That is,
\begin{align}
\nu_-^*:=-min\lbrace\Re(\mu_i(\lambda_{br})),\ i=1,2,3\rbrace.\label{EQ:ideal_weight_definition}
\end{align}

\begin{figure}[t]
\centering
\scalebox{1}{\subfloat{\label{fig:allspecta}\includegraphics{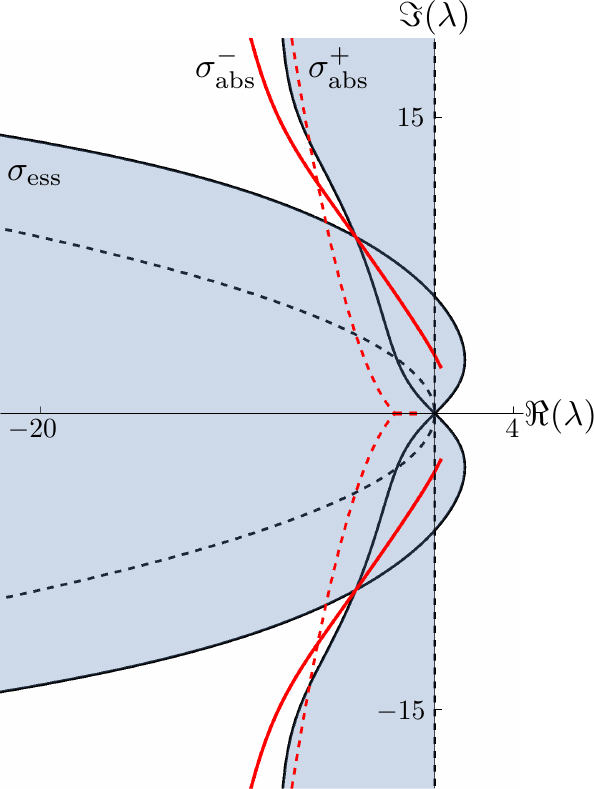}}}\hspace{.5cm}
\scalebox{1}{\subfloat{\label{fig:allspectb}\includegraphics{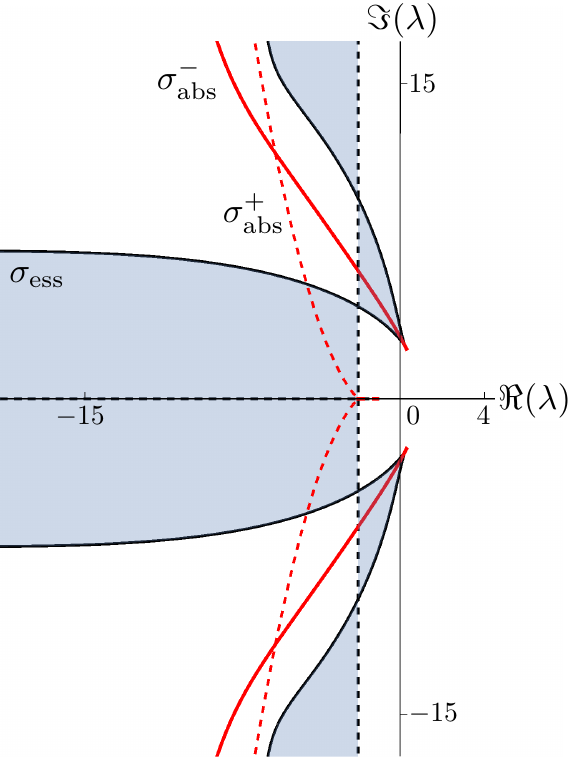}}}
\caption{The essential and absolute spectrum in the unweighted space (left panel) and in the ideally weighted space (right panel) for $\beta=c=2$, $\e=0$ and $m=0$, where the ideal weight is $\nu_-^*\approx -0.73$ and $\nu_+^*=c/2=1$. The dispersion relations of $M_0^++\nu_+I$ \eqref{EQ:dispersion_plus_weighted} are shown as black dashed lines, while those of $M_0^-+\nu_-I$ \eqref{EQ:dispersion_weighted_minus} are shown as black solid lines, $\sigma_{\rm abs}^+$ is shown as red dashed lines and $\sigma_{\rm abs}^-$ as red solid lines. The shaded regions are the interior of the (weighted) essential spectrum. Note the ideally weighted essential spectrum still contains values in the right half plane and the travelling wave solutions are thus absolutely unstable.}
\label{FIG:allspect}
\end{figure}

\begin{figure}[t]
\centering
\scalebox{1}{\subfloat{\label{fig:allspect_stablea1}\includegraphics{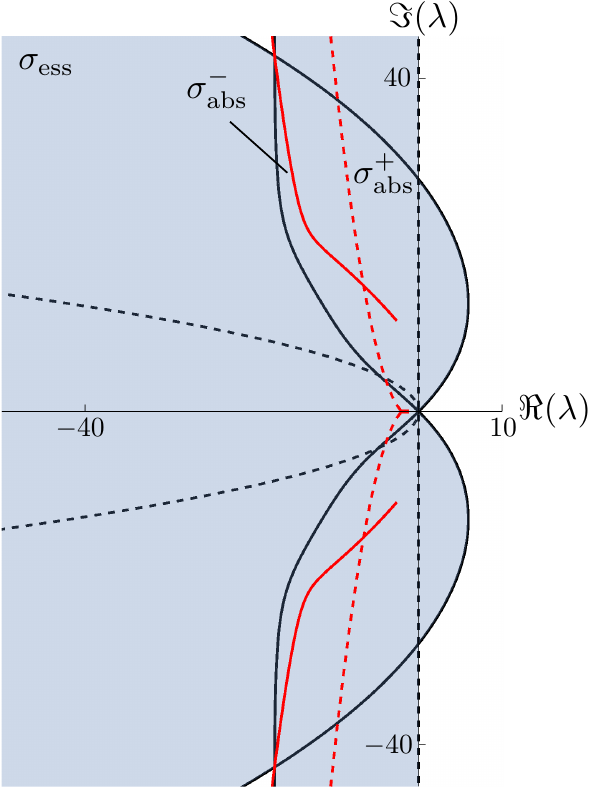}}}\hspace{.5cm}
\scalebox{1}{\subfloat{\label{fig:allspect_stableb2}\includegraphics{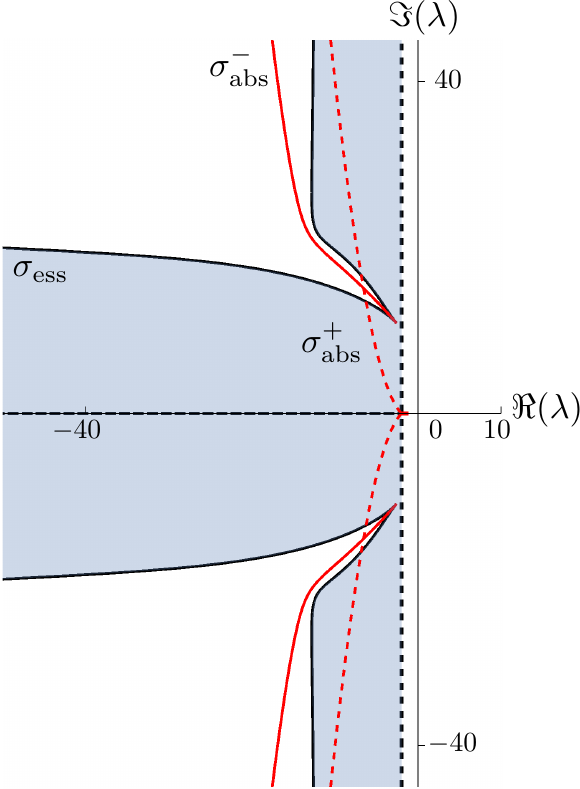}}}
\caption{
The essential and absolute spectrum in the unweighted space (left panel) and in the ideally weighted space (right panel) for 
$\beta=1.3<\beta_{\rm crit}$ \eqref{EQ:10th_order_poly}, $c=2$, $\e=0$ and $m=0$, where the ideal weight is $\nu_-^*\approx -2.445$ and $\nu_+^*=c/2=1$.
The dispersion relations of $M_0^++\nu_+I$ 
\eqref{EQ:dispersion_plus_weighted} are shown as black dashed lines, while those of $M_0^-+\nu_-I$ \eqref{EQ:dispersion_weighted_minus}
are shown as black solid lines, $\sigma_{\rm abs}^+$ is shown as red dashed lines and $\sigma_{\rm abs}^-$ as red solid lines. The shaded regions are the interior of the (weighted) essential spectrum. Note the ideally weighted essential spectrum is fully contained in the left half plane.}
\label{FIG:allspect_stable}
\end{figure}

We have outlined how to locate the full essential and absolute spectrum, as well as how to compute the ideal weights, for a given parameter set. See, for example, Figures \ref{FIG:allspect} and \ref{FIG:allspect_stable}. For the parameter values used in Figure \ref{FIG:allspect}, the ideally weighted essential spectrum and absolute spectrum contain values in the right half plane and the travelling wave solution is thus absolutely unstable. In contrast, for the parameter values used in Figure \ref{FIG:allspect_stable}, there exists a range of weights such that the essential spectrum (in the weighted space) is in the open left half plane and the travelling wave solution is potentially only transiently unstable.
Observe that $M_0^+$ requires positive weights $\nu_+$ to weigh its dispersion relations into the open left half plane, while $M_0^-$ requires negative weights $\nu_-$, necessitating the two-sided weight \eqref{EQ:two_sided_weight}.

 \begin{remark}\label{REM:branch_points}
 We refer to the value $\lambda$ such that $\mu(\lambda)$ is a second order root of \eqref{EQ:char_poly_minus} and $\lambda\in\sigma_{\rm abs}$ as a {\it branch point} because it is a branch point of the Evans function, an analytic tool used to locate the point spectrum. In general, not all spatial eigenvalues with algebraic multiplicity greater than one are contained in the absolute spectrum, they also occur in the generalised absolute spectrum. It is also not always the case that the leading edge of the absolute spectrum is a branch point, see for example \cite{sandstede2000absolute}. However, for the Keller-Segel model the leading edge of the sets $\sigma_{\rm abs}^\pm$ do coincide with branch points. See \S\ref{SS:POINT} for further discussion of the Evans function and point spectrum.
 \end{remark}

\subsection{Proof of Theorem~\ref{TH:MAIN1} for $\e=m=0$}
\label{SS:PROOF1}
From Figures \ref{FIG:allspect} and \ref{FIG:allspect_stable} it is clear that  
there is a transition from absolute spectrum fully contained in the left half plane to absolute spectrum entering into the right half plane. 
Consequently, there must be a critical set of parameters such that the branch point $\lambda_{br}$ solving \eqref{EQ:discriminant} is purely imaginary. So, we set $\lambda_{br}:= i\lambda$, $\lambda\in\mathbb{R}$, and equate the real and imaginary parts of \eqref{EQ:discriminant} to zero. This gives
\begin{align*}
\lambda ^4-\frac{\beta  \left(5 \beta ^3-28 \beta
   ^2+50 \beta -26\right) c^4 \lambda^2}{(\beta -1)^2 (2 \beta -1)^2}+\frac{\beta ^2 c^8}{(\beta -1)^2 (2 \beta -1)^2}&=0,\numberthis\label{EQ:br_cond1} \\
\lambda\left(\lambda^4-\frac{\beta  \left(18
   \beta ^2-37 \beta +20\right) c^4 \lambda ^2}{2 (\beta -1)^3} -\frac{\beta  \left(\beta ^2-6 \beta +2\right) c^8}{2 (\beta -1)^4}\right)&=0. \numberthis\label{EQ:br_cond2}
\end{align*}
Since $\lambda=0$ is not a solution of \eqref{EQ:br_cond1}, the transition occurs away from the real axis, \ie\ the branch points form a complex \edit{conjugate} pair.
Moreover, we can divide out $\lambda$ from \eqref{EQ:br_cond2} and the roots of \eqref{EQ:br_cond2} are given by $\lambda = \pm \sqrt{\Lambda_{1,2}}$ with
\begin{align}
\Lambda_{1,2}&=\frac{c^4\left(\beta  \left(18 \beta ^2-37 \beta +20\right)\pm \sqrt{\Delta} \right)}{4
(\beta -1)^3}, \label{eqn_L2}
\end{align}
where
\begin{align*}
\Delta&:=\beta 
\left(324 \beta ^5-1324 \beta ^4+2025 \beta ^3-1360 \beta ^2+320 \beta +16\right).
\end{align*} 
It follows from Sturm's Theorem, see, for instance, Theorem 6.3d in \cite{henrici1978applied}, that $\Delta>0$ for all $\beta>1$, \ie\ $\Lambda_{1,2}$ are real-valued for $\beta>1$.
Substituting these roots into \eqref{EQ:br_cond1} gives 
   \begin{align*}
&\frac{\beta  c^8}{8 (\beta -1)^4}\Big(1116 \beta ^7-5  050 \beta ^6+8422 \beta ^5-5440 \beta ^4-455 \beta ^3  \\
  &\left.+2104 \beta ^2 -704 \beta+8\pm\left(62
   \beta ^4-154 \beta ^3+90 \beta ^2+35 \beta -32\right) \sqrt{\Delta}\right)=0.
\end{align*}
Since $\beta>1$ and $c>0$, this is equivalent to
\begin{equation}
 \begin{split}
&\left(1116 \beta ^7-5  050 \beta ^6+8422 \beta ^5-5440 \beta ^4-455 \beta ^3 +2104 \beta ^2 -704 \beta+8\right)\\
&=\pm\left(-62
   \beta ^4+154 \beta ^3-90 \beta ^2-35 \beta +32\right)\sqrt{ \Delta},\label{EQ:beta_crit_before_squaring}
\end{split}
\end{equation}
which is independent of $c$, see Remark \ref{REM:scaling_invariance_c}. Squaring \eqref{EQ:beta_crit_before_squaring} gives
\begin{align*}
 16(\beta-1)^3f(\beta)&=16(\beta-1)^3\left(310 \beta ^{10}-3234 \beta ^9+17112 \beta ^8-49101 \beta ^7+76180 \beta ^6\right.\\
 &\quad\left.-58398 \beta
   ^5+10056 \beta ^4+15040 \beta ^3-9680 \beta ^2+1716 \beta -4\right)\\
   &=0,
   \end{align*}
   where $f(\beta)$ is the same polynomial as the polynomial \eqref{EQ:10th_order_poly} of Theorem \ref{TH:MAIN1}. 
So, the purely imaginary branch points indicating the transition to absolute instability are determined by the root $\beta_{\rm crit}$.
In particular, $\beta_{\rm crit}\approx1.6195$ solves \eqref{EQ:br_cond1} and \eqref{EQ:br_cond2} with $\lambda_{br}^\pm= \pm i \sqrt{\Lambda_1(\beta_{\rm crit})}\approx\pm 1.0883\,c^2\,i$. 
   
As there is only one root of \eqref{EQ:10th_order_poly} satisfying the condition $\beta>1$, the absolute spectrum is fully contained in the open left half plane for $1<\beta<\beta_{\rm crit}$, \ie\ the transition into the right half plane only happens for $\beta=\beta_{\rm crit}$.
Since the absolute spectrum always contains values in the right half plane for $\beta>\beta_{\rm crit}$, all travelling wave solutions with $\beta>\beta_{\rm crit}$ are absolutely unstable. This concludes the proof of Theorem~\ref{TH:MAIN1} for $\e=m=0$.

\begin{remark}\label{REM:weighted_dispersion_range}
It is possible for the absolute spectrum of an operator to be contained in the open left half plane, yet the weighted essential spectrum contains values in the right half plane for all weights. This is referred to as an {\it essential instability}, see \cite{sandstede2000absolute} for examples of essential instabilities. 
We now show that for a range of weights, the weighted dispersion relations, and thus the weighted essential spectrum, do not cross into the right half plane for $1< \beta<\beta_{\rm crit}$, \ie\ travelling wave solutions in the Keller-Segel model do not exhibit essential instabilities. The ideally weighted dispersion relations \eqref{EQ:dispersion_plus_max_weighted} and absolute spectrum $\sigma_{\rm abs}^+$ \eqref{EQ:abs_plus_m0} associated with $M_0^+$ are contained in the open left half plane for $1< \beta<\beta_{\rm crit}$. So, what remains to prove is that there exists a range of weights such that the weighted dispersion relations of $M_0^-$ \eqref{EQ:dispersion_weighted_minus} are fully contained in the open left half plane for $1< \beta<\beta_{\rm crit}$.

The characteristic polynomial of $M_0^-+\nu_- I$ \eqref{EQ:dispersion_weighted_minus} is quadratic in $\lambda\in\mathbb{C}$.
So, we can explicitly solve for $\lambda_{1,2}$ and extract the real parts of the solutions.
It follows that
 \begin{align}
 \label{EQ:limit22}
 &\lim_{|k|\to\infty} \Re(\lambda_1)=-c\left(\frac{c\beta  }{\beta -1}+\nu_-\right),\qquad  \lim_{|k|\to\infty} \Re(\lambda_2)=-\infty.
 \end{align}
That is, the dispersion relations of $M_0^-+\nu_- I$ approach vertical lines in the complex plane. Requiring that $\Re(\lambda_1)<0$ as $|k|\to\infty$ gives a lower bound on admissible weights $\nu_->-\frac{c\beta}{\beta-1}$ (note that it turns out that this lower bound is not sharp, see Figure~\ref{FIG:range_permissable}).

Next, we compute the values $\lambda$ where the dispersion relations of $M_0^-+\nu_- I$ \eqref{EQ:dispersion_weighted_minus} cross the imaginary axis. 
Therefore, we assume that $\lambda$ is purely imaginary and solve \eqref{EQ:dispersion_weighted_minus}. 
This way, we eliminate the parameter $k$ and obtain a 
cubic polynomial equation in $\Lambda:=\Im(\lambda)^2$ (with unknowns $\beta, c$ and $\nu_-$). 
So, non-negative real roots of this polynomial in $\Lambda$ correspond to the intersections of the dispersion relations of $M_0^-+\nu_- I$ with the imaginary axis.
In the unweighted case $\nu_-=0$ it has one positive root and a root in the origin, see also the left panels of Figure~\ref{FIG:allspect} and \ref{FIG:allspect_stable}.   
For decreasing $\nu_-$, these two roots approach each other and collide at $\nu_{\rm max}=\nu_{\rm max}(\beta, c)$ (while the third root stays negative). 
The polynomial has no non-negative real roots if we further decrease $\nu_-$. 
These weights correspond to the case where the weighted dispersion relations do not intersect the imaginary axis and are thus fully contained in the open left half plane. 
At $\nu_{\rm min}=\nu_{\rm min}(\beta, c)$ two positive roots reappear (while the third root is still negative) and these positive roots persist upon further decreasing $\nu_-$. 
In other words, for weights $\nu_- \in(\nu_{\rm min},\nu_{\rm max})$ the dispersion relations of $M_0^-+\nu_- I$ \eqref{EQ:dispersion_weighted_minus} never intersect the imaginary axis and are fully contained in the open left half plane.
The values $\nu_{\rm min}$ and $\nu_{\rm max}$ are given as the roots of an $11^{th}$ order polynomial in $\nu_-$ and the range of admissible weights shrinks to a point as $\beta\uparrow\beta_{\rm crit}$, see Figure \ref{FIG:range_permissable}. In particular, one rediscovers $f$ \eqref{EQ:10th_order_poly} by equating the derivative of this $11^{th}$ order polynomial to zero. This is equivalent to finding $\beta$ such that $\nu_{\rm min}=\nu_{\rm max}$. 
Obtaining the range of admissible weights is straightforward for given values of $\beta$ and $c$, but complicated to determine for general $1<\beta<\beta_{\rm crit}$ and $c$. 
See Figure \ref{FIG:range_permissable} for a plot of $\nu_{\rm max}$ and $\nu_{\rm min}$ (and the ideal weight $\nu_-^*$ obtained from \eqref{EQ:ideal_weight_definition}) versus $\beta$.

 \begin{figure}[t]
\centering
\scalebox{1.1}{\includegraphics{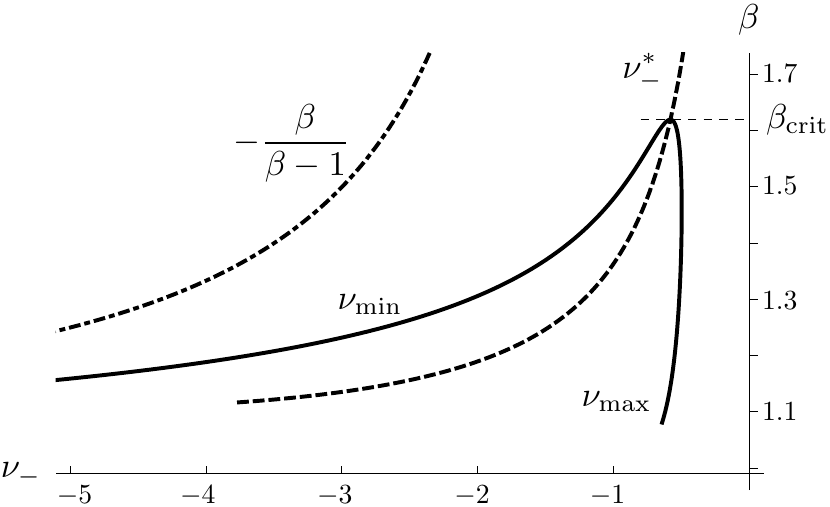}}
\caption{The area in between the curves $\nu_{\rm min}$ and $\nu_{\rm max}$ indicates the range of weights $\nu_-$ such that the absolute spectrum and weighted essential spectrum are contained in the open left half plane for $c=1,\ m=0$, $\e=0$ (and $\nu^*_+ = c/2=1/2$). 
For $\beta=\beta_{\rm crit}$ the values $\nu_{\rm min},\ \nu_{\rm max}$ and the ideal weight $\nu_-^*$ coincide and so the essential spectrum cannot be weighted into the open left half plane for $\beta \geq \beta_{\rm crit}$. 
The dot-dashed curve represents the asymptotic condition $\nu_->-\frac{\beta}{\beta-1}$ coming from \eqref{EQ:limit22}. }
\label{FIG:range_permissable}
\end{figure}
\end{remark}

\begin{remark}\label{REM:scaling_invariance_c}
The results on the existence of a range of weights to move the essential spectrum into the open left half plane and the (in)stability of the absolute spectrum are independent of the wave speed $c$. This is not a coincidence as the dispersion relations can be rescaled to be independent of $c$. In particular, the substitutions $\lambda=c^2\tilde{ \lambda},\ \nu=c\tilde{ \nu},\  k=c\tilde{k}$ transform the dispersion relations of $M_0^++\nu_+ I$ \eqref{EQ:dispersion_plus_weighted} into
\begin{equation*}
c^2\tilde{\lambda}=c^2\left(-\tilde{\nu}_++i\tilde{k}\right),\ \text{ and }\ c^2\tilde{\lambda}=c^2\left(-\tilde{k}^2-\tilde{\nu}_+(1-\tilde{\nu})+i(\tilde{k}-2\tilde{k}\tilde{\nu}_+\right),
\end{equation*}
which is equivalent to the dispersion relations of $M_0^++\nu_+ I$ \eqref{EQ:dispersion_plus_weighted} for $c=1$.
Similarly, the dispersion relations of $M_0^-+\nu_- I$ \eqref{EQ:dispersion_weighted_minus} become
  \begin{equation*}
 \begin{aligned}
 &c^4\bigg(\tilde{\lambda}^2+\left(\frac{(2-\beta ) (i \tilde{k}-\tilde{\nu}_-)}{\beta -1}-(i \tilde{k}-\tilde{\nu}_-)^2\right)\tilde{\lambda}+\frac{\beta   (i \tilde{k}-\tilde{\nu}_-)}{(\beta -1)^2}\\
 &\qquad-\frac{(\beta +1)  (i \tilde{k}-\tilde{\nu}_-)^2}{\beta -1}+ (i \tilde{k}-\tilde{\nu}_-)^3\bigg)=0,
 \end{aligned}
 \end{equation*}
 which is equivalent to the dispersion relations of $M_0^-+\nu_- I$ \eqref{EQ:dispersion_weighted_minus} for $c=1$.
 In other words, the magnitude of $c$ does not affect the (in)stability results and only affects the multiplicative scaling of the spectrum. As a consequence, all the figures presented in this manuscript are {\it generic} in $c$ up to the above scaling of $\lambda, \nu$ and $k$.
 \end{remark}

\section{Sublinear and linear consumption and zero diffusivity of the chemoattractant}
\label{S:SUBLIN}
In this section, we examine the effect of the parameter $m$ on the location of the weighted essential spectrum and absolute spectrum associated with a travelling wave solution. 
Since travelling wave solutions only exist for $0\leq m\leq 1$, {\it e.g.} \cite{wang2013mathematics}, we take $0<m\leq1$. 
We prove Theorem~\ref{TH:MAIN1} for $0<m\leq1$ and $\e=0$.
It turns out that the analysis for $0< m<1$ is similar, at least qualitatively, to the analysis of the previous section for $m=0$. The analysis simplifies significantly for $m=1$ and \edit{we} note that the results of this case can be in part deduced from \cite{meyries2011local} where a version of the Keller-Segel model with nonzero growth rate is studied.
 
In particular, we show that for sublinear consumption, \ie\ $0<m<1$, there exists a critical value $\beta_{\rm crit}^m=\beta_{\rm crit}(1-m)$ (with $\beta_{\rm crit}$ the root of \eqref{EQ:10th_order_poly}) such that for $1-m<\beta<\beta_{\rm crit}^m$ the absolute spectrum is fully contained in the open left half plane.
The absolute spectrum enters the right half plane for
  $\beta>\beta_{\rm crit}^m$ and all travelling wave solutions are thus absolutely unstable for
  $\beta>\beta_{\rm crit}^m$.
For linear consumption, \ie\ $m=1$, we show that the absolute spectrum always contains the origin. 
Consequently, the essential spectrum cannot be weighted into the open left half plane.
\subsection{Set-up}
For $0<m\leq 1$ and $\e=0$, the eigenvalue problem is given by \eqref{EQ:eigen}, which we restate for convenience 
 \begin{align}&\mL \begin{pmatrix}
 p\\q
 \end{pmatrix}=\lambda\begin{pmatrix}
 p\\q
 \end{pmatrix},\quad\quad
\mathcal{L}:=\begin{pmatrix}c\pder{}{z}{}- m w u^{m-1}& -u^m\\ \mathcal{L}_p &\mathcal{L}_q
\end{pmatrix}\label{EQ:Lop_sublinear_m}
\end{align}
with
\begin{equation}
\label{EQ:Lop_sublinear_m_Lp_Lq}
\begin{aligned}
\mathcal{L}_p&:=\beta\left(\frac{w_z u_z }{u^2}+\frac{wu_{zz} }{u^2}-\frac{2wu_{z}^2}{u^3}\right)+\beta\left(\frac{2wu_{z}}{u^2}-\frac{w_z}{u}\right)\pder{}{z}{}-\frac{\beta w}{u}\pder{}{z}{2},\\
\mathcal{L}_q&:=\beta\left(\frac{u_{z}^2}{u^2}-\frac{u_{zz}}{u}\right)+\left(c-\frac{\beta u_{z}}{u}\right)\pder{}{z}{}+\pder{}{z}{2},
\end{aligned}
\end{equation}
where $u$ and $w$ are the travelling wave solutions given in \eqref{tw_profiles_0}. Observe that the first row of $\mL$ simplifies significantly in the cases $m=0$ and $m=1$.
We take a slightly different approach as in \S\ref{S:MZERO} and first write \eqref{EQ:Lop_sublinear_m} as a third order equation in $p$, see Remark \ref{REM:pq_subs}. From the first row of \eqref{EQ:Lop_sublinear_m} we have
\begin{equation}
q= c u^{-m}p_z - (mwu^{-1}+ \lambda u^{-m}) p\,,\label{EQ:q_m_nonzero}
\end{equation}
and we differentiate this to obtain
\begin{equation}
\label{EQ:qz_m_nonzero}
\begin{aligned}
q_z&=c u^{-m}p_{zz} +((c u^{-m})_z-(mwu^{-1}+ \lambda u^{-m}))p_z\\
&-\left(mwu^{-1}+ \lambda u^{-m}\right)_z p,\\
q_{zz}&=c u^{-m}p_{zzz}+\left(2(c u^{-m})_z-(mwu^{-1}+ \lambda u^{-m})\right)p_{zz}\\
&+\left((cmu^{-m})_{zz}-2(mwu^{-1}+ \lambda u^{-m})_z\right)p_z-(mwu^{-1}+ \lambda u^{-m})_{zz} p.
 \end{aligned}
 \end{equation}
We substitute \eqref{EQ:q_m_nonzero} and \eqref{EQ:qz_m_nonzero} into the second row of \eqref{EQ:Lop_sublinear_m}, that is into 
$ \mL_pp+\mL_qq=\lambda q$, and we eliminate $w$ using $w=cu_z u^{-m}$ (\eqref{EQ:KStw2} with $\e=0$). The resulting third order operator is
\begin{align}
p_{zzz}- \cC_mp_{zz}- \cB_mp_z - \cA_mp = 0 \label{EQ:Pop_m_nonzero_eps_zero}
\end{align}
where 
\begin{equation}
\begin{split}
\cA_m =  \, \, & \left(\lambda  (m+1) (\beta +m)-c^2 m\right)\frac{u_z^2 }{c u^2}+2 (m+1) (\beta
   +m) \frac{u_z^3}{u^3}-2\lambda m\frac{ u_z}{u}\\
   &-(2 \beta+3 m)\frac{ u_z u_{zz}}{u^2}-\lambda(\beta +m)\frac{  u_{zz}}{c u}-\frac{\lambda ^2}{c},\\
\cB_m =  \, \,  & \left(2 c^2 m-\lambda  (\beta +2 m)\right)\frac{u_z }{c u} -3 (m+1) (\beta +m) \frac{u_z^2}{u^2}+(2 \beta +3 m)\frac{ u_{zz}}{u}+2 \lambda,  \\
\cC_m =  \, \, & \frac{\lambda }{c}-c+(2 \beta +3 m)\frac{ u_z}{u}. 
\end{split}
\end{equation}
 Next, we set $p_1:=p_z$ and $p_2:=p_{zz}$, and define the operator $\mathcal{T}_m(\lambda)$ 
 \begin{equation}
\begin{split}
\mathcal{T}_m(\lambda)\begin{pmatrix}
p\\p_1\\p_2
\end{pmatrix} &:= \begin{pmatrix}
p\\p_1\\p_2
\end{pmatrix}' -M_m(z,\lambda) \begin{pmatrix}
p\\p_1\\p_2
\end{pmatrix}=0, \\  \\ M_m(z,\lambda) &:=\begin{pmatrix}
0&1&0\\ 0&0&1\\ \cA_m&\cB_m&\cC_m
\end{pmatrix}.\label{EQ:Top_m_nonzero}
\end{split}
\end{equation}
While we have used a slightly different approach compared to \S\ref{S:MZERO}, the spectrum of $\mathcal{T}_0(\lambda)$ in \eqref{EQ:T0} and the spectrum of $\mathcal{T}_m(\lambda)$ \eqref{EQ:Top_m_nonzero} agree in the limit $m\to0$. 

\begin{remark}\label{REM:pq_subs}
 The substitutions \eqref{EQ:q_m_nonzero} and \eqref{EQ:qz_m_nonzero} are necessary due to the appearence of the term $w/u$ appearing in $\mL_p$ \eqref{EQ:Lop_sublinear_m_Lp_Lq}. While the term $w/u$ is bounded for $m=0$, the term is unbounded as $z\to-\infty$ for $0<m<1$. 
 However,  
by making the substitutions \eqref{EQ:q_m_nonzero} and \eqref{EQ:qz_m_nonzero}, we obtain   \eqref{EQ:Pop_m_nonzero_eps_zero}, which is {\it asymptotically constant} and equivalent to $\mL$ \eqref{EQ:Lop_sublinear_m}. The equivalence of \eqref{EQ:Pop_m_nonzero_eps_zero} and $\mL$ \eqref{EQ:Lop_sublinear_m} becomes clearer when we see that \eqref{EQ:Pop_m_nonzero_eps_zero} is actually the linearised eigenvalue problem
obtained from eliminating $w(z,t)=u^{-m}(z,t)\left( c u_z(z,t)-u_t(z,t)\right)$ from \eqref{EQ:KStw} first.
\end{remark}
\subsection{Essential spectrum}
\label{SS:ESS_SPECT_MNONZERO}
\begin{figure}[t]
  \subfloat{\scalebox{1}{\includegraphics{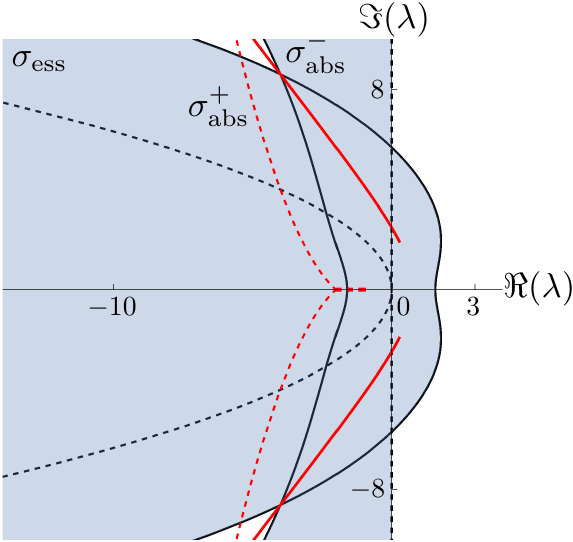}}}\hspace{.5cm}
\subfloat{\scalebox{1}{\includegraphics{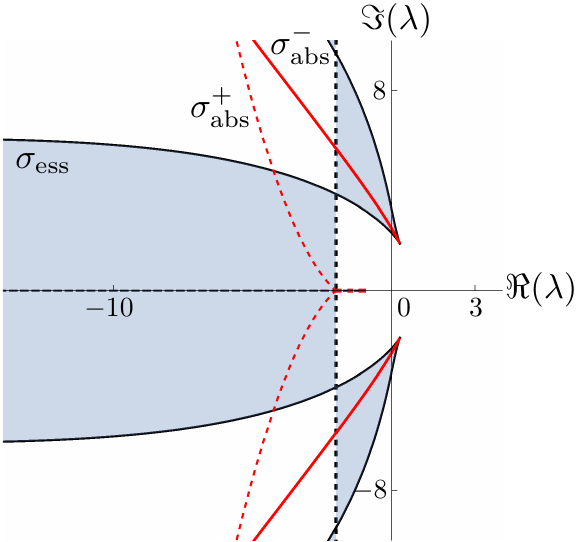}}}\\
\subfloat{\scalebox{1}{\includegraphics{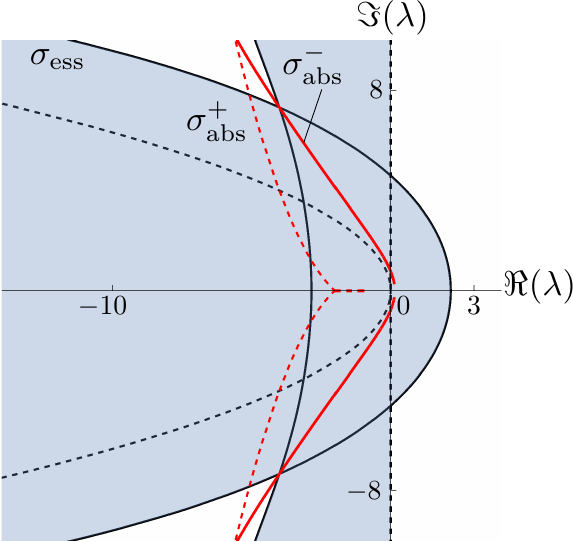}}}\hspace{.5cm}
\subfloat{\scalebox{1}{\includegraphics{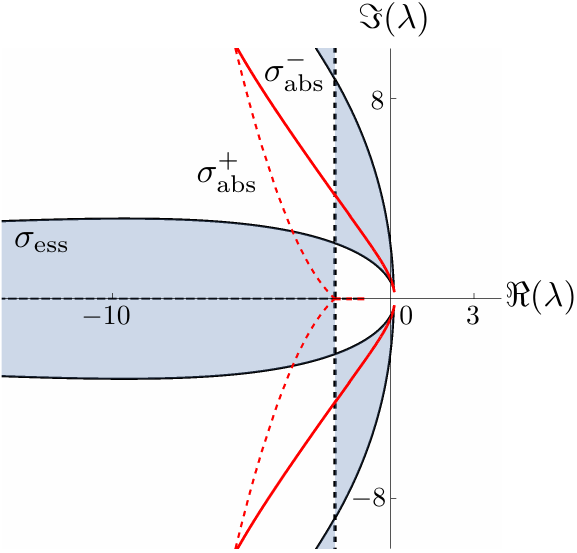}}}
  \caption{The essential and absolute spectrum for $m=0.1$ (upper panels) and $m=0.7$ (lower panels)
with  
  $\beta=c=2$ and $\e=0$. The dispersion relations of $M_m^{+}+\nu_+I$ (dashed black) and $\sigma_{\rm abs}^{m,+}$ (dashed red) are the same in all four panels and the ideal weight for $z\to\infty$ is still given by $\nu_+^*=c/2=1$. The dispersion relations of $M_m^{-}+\nu_-I$ are shown as solid black lines and $\sigma_{\rm abs}^{m,-}$ as solid red. 
  Upper left panel: the spectrum in the unweighted space for $m=0.1$. 
  Upper right panel: the ideally weighted space for $m=0.1$, where the ideal weight is $\nu_-^*\approx-0.778$. 
  Lower left panel:  the spectrum in the unweighted space for $m=0.7$. 
  Lower right panel: the ideally weighted space for $m=0.7$, where $\nu_-^*\approx-0.959$. 
  As $m$ increases to one, the real and imaginary components of the branch points  $\lambda_{br}^\pm$ decrease and approach the origin, see \S\ref{SS:MLINEAR}. }
\label{FIG:m_changing}
\end{figure}
We use the limits given in \eqref{EQ:limits} (with $\e=0$) and the fact that $u_{zz}=(w u^{m})_z/c$ \eqref{EQ:KStw2}, to 
obtain 
\begin{align*}
\lim_{z\to-\infty}\frac{u_z}{u}=\frac{c}{\beta+m-1},\quad \lim_{z\to-\infty}\frac{u_{zz}}{u}=\frac{c^2}{(\beta+m-1)^2},\quad \lim_{z\to\infty} (u,u_z,u_{zz})=(1,0,0).
\end{align*}
Using these limits, the asymptotic values of $\cA_m,\ \cB_m$ and $\cC_m$ as $z\to\pm\infty$, denoted $\cA_m^{\pm},\ \cB_m^{\pm}$ and $\cC_m^\pm$ respectively, are
\begin{align}
 \cA_m^{+}= \, \,& -\frac{\lambda ^2}{c},\quad \cB_m^{+}=    \, \,  2\lambda, \quad \cC_m^{+}=   \, \,  \frac{\lambda }{c}-c,
\end{align}
and
\begin{equation}
\begin{split}
 \cA_m^{-}= \, \, &-\frac{\lambda ^2}{c}-\frac{c   m (\beta +m-2)}{(\beta +m-1)^2}\lambda+\frac{c^3 m (\beta +m)}{(\beta +m-1)^3},\\ 
 \cB_m^{-}=    \, \, & \frac{(\beta -2)  }{\beta +m-1}\lambda-\frac{c^2 \left(\beta +m(\beta+m  +2) \right)}{(\beta +m-1)^2}, \\ 
 \cC_m^{-}=   \, \, & \frac{\lambda }{c}+\frac{c (\beta +2 m+1)}{\beta +m-1}.
\end{split}
\end{equation}
We define the asymptotic matrices
\begin{align}
\label{EQ:M}
 M_m^{\pm}(\lambda):=\displaystyle\lim_{z\rightarrow\pm\infty}M_m(z,\lambda)=\begin{pmatrix}
 0 & 1 & 0 \\
 0 & 0 & 1 \\
 \mathcal{A}_m^{\pm} &  \mathcal{B}_m^{\pm} &  \mathcal{C}_m^{\pm}
\end{pmatrix},
\end{align}
related to the asymptotic operator associated with $\mathcal{T}_m$ \eqref{EQ:Top_m_nonzero}.
The dispersion relations of $M_m^{+}$ are independent of $m$ and $\beta$, and the same as for $m=0$ \eqref{EQ:dispersion_plus}. 
The dispersion relations of $M_m^{-}$ depend on $m$ and are implicitly given by
\begin{equation}
\begin{split}
&\lambda^2+\left(k^2+\frac{c^2 m (\beta +m-2)}{(\beta +m-1)^2}-\frac{ick (\beta -2)}{\beta +m-1}\right)\lambda+\frac{c^2 k^2 (\beta +2
   m+1)}{\beta +m-1}\\
&-\frac{c^4 m (\beta +m)}{(\beta +m-1)^3}+ \frac{ic^3k  \left(\beta +m(m+\beta +2)\right)}{(\beta +m-1)^2}-ick^3=0.\label{EQ:dispersion_m_nonzero_minus}
\end{split} 
\end{equation}
In the limit $m\to0$, \eqref{EQ:dispersion_m_nonzero_minus} coincides with the dispersion relations of $M_0^-$   \eqref{EQ:dispersion_minus}. The dispersion relations $M_m^{+}$ \eqref{EQ:dispersion_plus} and $M_m^{-}$ \eqref{EQ:dispersion_m_nonzero_minus} form the boundaries of the essential spectrum and 
$\lambda\in\mathbb{C}$ such that \ie\ $i_+\neq i_-$ (see Definition \ref{DEFN:EssSpect}) forms
the interior of the (unweighted) essential spectrum. 
See the two left panels of Figure~\ref{FIG:m_changing} for the unweighted essential spectrum for two different values of $m$.
\subsection{The weighted essential spectrum and the absolute spectrum}

As for $m=0$, we consider a two-sided weight of the form \eqref{EQ:two_sided_weight}.
Since the dispersion relations of $M_m^{+}$ and $M_0^+$ are the same
the ideal weight for $z \to \infty$ are unchanged for $0<m \leq 1$. That is, $\nu_+^*=c/2$. Consequently, $\sigma_{\rm abs}^{m,+}=\sigma_{\rm abs}^+$ \eqref{EQ:abs_plus_m0}. See also Figure \ref{fig:absplus}.

The dispersion relations of $M_m^{-}+\nu_- I$ are implicitly given by
\begin{equation}
\begin{split}
&\lambda ^2+ \lambda  \left(-(ik-\nu_-)^2+\frac{c^2 m (\beta +m-2)}{(\beta +m-1)^2}-\frac{c(i
   k-\nu_-)(\beta-2 ) }{\beta +m-1}\right)\\
   &-\frac{c^2 (ik-\nu_-)^2 (\beta +2
   m+1)}{\beta +m-1}-\frac{c^4 m (\beta +m)}{(\beta +m-1)^3}\\
   &+\frac{c^3 (ik-\nu_-) (\beta +m (\beta +m+2))}{(\beta +m-1)^2}+c(ik-\nu_-)^3=0.\label{EQ:dispersion_m_nonzero_minus_weighted}
\end{split}
\end{equation}
The shift in the essential spectrum due to weighting in the $0<m\leq1$ case is qualitatively similar to the behaviour shown in Figure \ref{FIG:various_weighted_minus}. That is, under a large range of weights the dispersion relations have self-intersections 
and these self-intersections form part of the absolute spectrum $\sigma_{\rm abs}^{m,-}$.  
So,
we can once again use a find root procedure on the weighted dispersion relations \eqref{EQ:dispersion_m_nonzero_minus_weighted} to locate $\sigma_{\rm abs}^{m,-}$.
See Figure~\ref{FIG:m_changing} for the unweighted essential spectrum, the ideally weighted essential spectrum, and the absolute spectrum for two different values of $m$.

\subsection{Proof of Theorem \ref{TH:MAIN1} for $0<m<1$ and $\e=0$}
\label{SS:PROOF3}
For $0<m< 1$ and $\e=0$, a polynomial $f_m(\beta)$, similar to the polynomial $f(\beta)$ \eqref{EQ:10th_order_poly} for $m=0$, can be derived. Its root $\beta_{\rm crit}^m = \beta_{\rm crit}(1-m)>1-m$ predicts the transition of the absolute spectrum into the right half plane (for increasing $\beta$).  
For $1-m<\beta<\beta_{\rm crit}^m$, the absolute spectrum is fully contained in the open left half plane. 
For $\beta>\beta_{\rm crit}^m$, the absolute spectrum enters the right half plane and the travelling wave solutions are thus absolutely unstable.

To determine the transition of the absolute spectrum into the right half plane we follow the same procedure as in \S\ref{SS:PROOF1} and we treat the characteristic polynomial of $M_m^{-}$ as a cubic polynomial in $\mu$ and equate the discriminant to zero. This gives
\begin{equation}
\label{EQ:disc_m}
\begin{aligned}
&\lambda^5+\frac{c^2 (2 \beta +m-1)^2}{4 (\beta +m-1)^2}\lambda^4+\frac{\beta  c^4 \left(18 \beta
   ^2+37 \beta  (m-1)+20 (m-1)^2\right)}{2 (\beta +m-1)^3}\lambda^3\\
   &+\frac{\beta  c^6 \left(5 \beta ^3+28 \beta ^2 (m-1)+50 \beta  (m-1)^2+26 (m-1)^3\right)}{4 (\beta +m-1)^4}\lambda^2\\
   &+\frac{\beta  c^8 (m-1) \left(\beta ^2+6 \beta  (m-1)+2 (m-1)^2\right)}{2 (\beta+m-1)^4}\lambda+\frac{\beta ^2 c^{10} (m-1)^2}{4 (\beta +m-1)^4}=0.
\end{aligned}
\end{equation}
This discriminant has a purely imaginary root under the condition
\begin{equation}
\label{EQ:beta_crit_m_nonzero}
\begin{aligned}
0&=\frac{\beta ^2 c^{20} (m-1)}{64 (\beta +m-1)^{13}}f_m(\beta)
\end{aligned}
\end{equation}
where 
\begin{equation}
\begin{aligned}
&f_m(\beta):=\left(310 \beta ^{10}+3234 \beta ^9 (m-1)+17112 \beta ^8 (m-1)^2+49101 \beta ^7 (m-1)^3\right.\\
&\quad+76180 \beta ^6 (m-1)^4+58398 \beta ^5 (m-1)^5+10056 \beta ^4
   (m-1)^6\\
   &\quad\left.-15040 \beta ^3 (m-1)^7-9680 \beta ^2 (m-1)^8-1716 \beta  (m-1)^9-4 (m-1)^{10}\right)\,.
\end{aligned}
\end{equation}
For $m=1$, \eqref{EQ:beta_crit_m_nonzero} is trivially satisfied. Therefore, we treat the $m=1$ case seperately, see \S\ref{SS:MLINEAR}. Upon introducing the variable $B=\frac{\beta}{(1-m)}$ (and setting $0<m<1$), \eqref{EQ:beta_crit_m_nonzero} becomes,
\begin{equation}
\begin{aligned}
0&=\frac{-B^2 c^{20}}{64 (B-1)^{13}}\left(310 B^{10}-3234 B^9+17112 B^8-49101 B^7+76180 B^6\right.\\
&\quad\left.-58398 B^5+10056 B^4+15040 B^3-9680 B^2+1716 B-4\right)\\
&=\frac{-B^2 c^{20}}{64 (B-1)^{13}}f(B),
\end{aligned}
\end{equation}
where $f$ is given by \eqref{EQ:10th_order_poly}. So, the roots of $f_m$ and $f$ are related by $\beta_{\rm crit}^m=\beta_{\rm crit}(1-m)$, and $\beta_{\rm crit}^m$ is the only root of \eqref{EQ:beta_crit_m_nonzero} that satisfies the condition $\beta+m>1$. In conclusion, we have that the absolute spectrum is fully contained in the open left half plane for $0\leq m<1, \e=0$ and $1-m<\beta<\beta_{\rm crit}^m$, while the absolute spectrum enters into the right half plane for $0\leq m<1, \e=0$ and $\beta>\beta_{\rm crit}^m$.  
This concludes the proof of Theorem \ref{TH:MAIN1} for $0<m<1$ and $\e=0$.

\begin{remark}
\label{REM:weighted_dispersion_range2}
Similar to the $m=0$ case, there also exist a range of weights $\nu_{\rm min}^m<\nu_-<\nu_{\rm max}^m$ for $0<m<1$ and $\e=0$,  such that the weighted essential spectrum is contained in the open left half plane for $1-m<\beta<\beta_{\rm crit}^m$. In other words, there are no essential instabilities in this case.
See also Remark \ref{REM:weighted_dispersion_range}.
\end{remark} 
\subsection{Linear consumption}
\label{SS:MLINEAR}
In the case of linear consumption, \ie\ $m=1$, the travelling wave solutions $(u,w)$ \eqref{tw_profiles_0} are a pair of wavefronts, rather than a pulse and a wavefront, see, for example, the right panel of Figure \ref{fig:tw_profiles}. 
In this case, the absolute spectrum and the ideally weighted essential spectrum contain the origin for all $\beta$ and as a result the essential spectrum cannot be weighted into the open left half plane. 
Note that the model \eqref{EQ:KStw1} with $m=1$ can be seen as a limit case of the model with non-zero growth term ($\kappa>0$ in \eqref{EQ:KSgen}) considered in \cite{meyries2011local}.

\begin{figure}[t]
  \centering
\scalebox{1}{\subfloat{\label{fig:aaa}\includegraphics{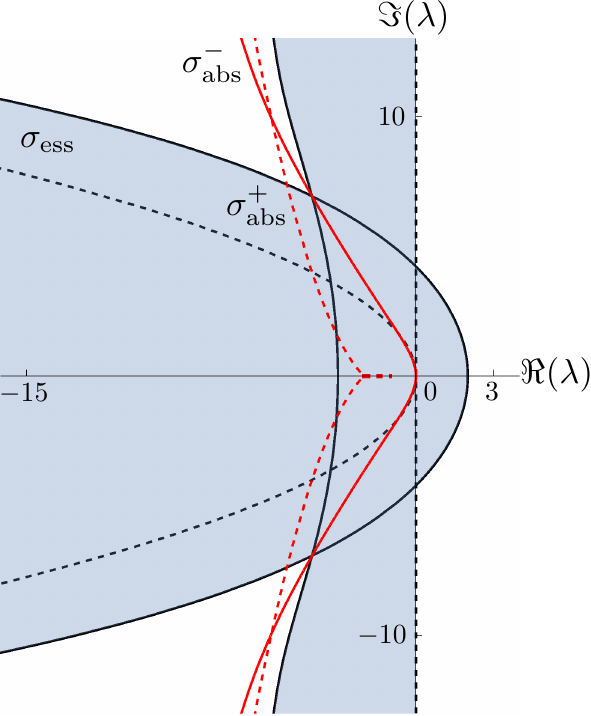}}}\hspace{.5cm}
\scalebox{1}{\subfloat{\label{fig:bbb}\includegraphics{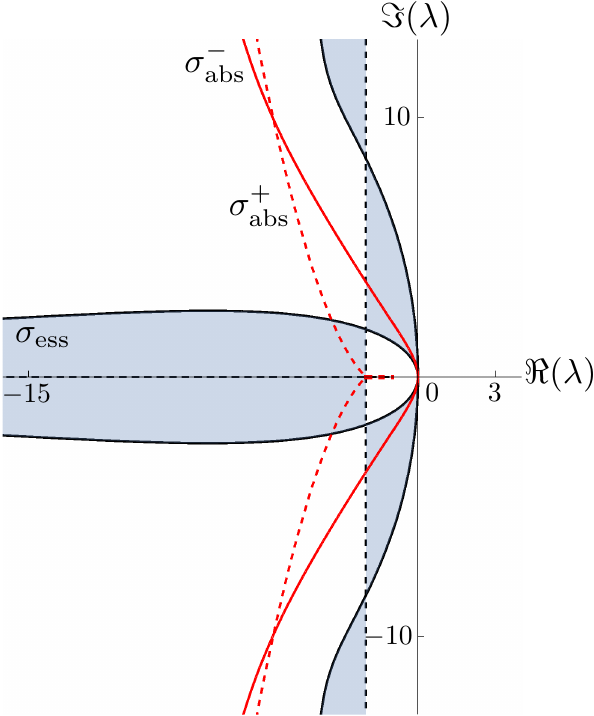}}}
  \caption{The essential and absolute spectrum in the unweighted space (left panel) and in the ideally weighted space (right panel) for $\beta=c=2$, $\e=0$ and $m=1$, where the ideal weight is $\nu_-^*=-c/\beta=-1$ and $\nu_+^*=c/2=1$. The dispersion relations of $M_1^{+}+\nu_+I$ \eqref{EQ:dispersion_plus_weighted} are shown as black dashed lines, while those of $M_1^{-}+\nu_-I$ \eqref{EQ:dispersion_m_nonzero_minus_weighted} are shown as black solid lines, $\sigma_{\rm abs}^{1,+}$ is shown as red dashed lines and $\sigma_{\rm abs}^{1,-}$ as red solid lines. The shaded regions are the interior of the (weighted) essential spectrum. The absolute spectrum contains the origin (for all parameter values $\beta$ and $c$) and the essential spectrum thus cannot be weighted into the open left half plane. }  \label{FIG:allspect_m_1}
\end{figure}

The dispersion relations of $M_1^{+}$ are independent of $m$ and $\beta$, see \S\ref{SS:ESS_SPECT_MNONZERO}, and therefore $\sigma_{\rm abs}^{1,+}=\sigma_{\rm abs}^+$ \eqref{EQ:abs_plus_m0} is fully contained in the open left half plane. Consequently, we only need to examine $\sigma_{\rm abs}^{1,-}$. 
The characteristic polynomial of $M_1^{-}$ is 
\begin{equation}
 \begin{aligned}
&\mu ^3-\mu ^2 \left(\frac{\beta\left(\beta +3  \right) c}{\beta
   ^2}+\frac{\lambda }{c}\right)+\mu  \left(\frac{(2-\beta ) \lambda }{\beta }+\frac{\beta\left(\beta ^2+(\beta -1)  +4 \beta \right) c^2}{\beta ^3}\right)\\
   &-\frac{(\beta +1) c^3}{\beta ^3}+\frac{(\beta -1) c \lambda }{\beta ^2}+\frac{\lambda ^2}{c}=0.\label{EQ:char_poly_m_1}
 \end{aligned}
 \end{equation}
To locate $\sigma_{\rm abs}^{1,-}$, we follow the same process as for $0\leq m<1$. In particular, we locate $\lambda\in\sigma_{\rm abs}^{1,-}$ such that the characteristic polynomial \eqref{EQ:char_poly_m_1} has a second order root in $\mu$. 
That is, we locate the branch points $\lambda_{br}^\pm$.
We equate the discriminant of \eqref{EQ:char_poly_m_1} to zero to obtain
 \begin{align}
 \lambda ^2 \left(4 \lambda ^3+4 c^2 \lambda ^2+36 c^4 \lambda +5 c^6 \right) = 0,\label{EQ:discriminant_m_1}
 \end{align} 
which has a second order root $\lambda=0$.
For $\lambda=0$, \eqref{EQ:char_poly_m_1} becomes 
\begin{align}
(\beta  \mu-c(\beta+1)) (c-\beta  \mu )^2=0 \quad \implies \quad \mu_{1}=\frac{(\beta+1)c}{\beta},\quad\mu_{2,3}=\frac{c}{\beta}.
\end{align}
Since $\Re(\mu_1)>\Re(\mu_2)=\Re(\mu_3)$, $0 \in \sigma_{\rm abs}^{1,-}$ and the ideal weight is $\nu_-^*=-\Re(\mu_{2,3})=-\frac{c}{\beta}$ \eqref{EQ:ideal_weight_definition}.
Furthermore, the ideally weighted essential spectrum and the absolute spectrum contain the origin for all $\beta$. That is, there are no parameter values such that the essential spectrum can be weighted fully into the open left half plane, see, for example, Figure \ref{FIG:allspect_m_1}. 
Note that the other three roots of \eqref{EQ:discriminant_m_1} are part of the generalised absolute spectrum.
This concludes the proof of Theorem \ref{TH:MAIN1} for $m=1$ and $\e=0$.

\begin{figure}[t]

  \centering
\scalebox{1}{\includegraphics{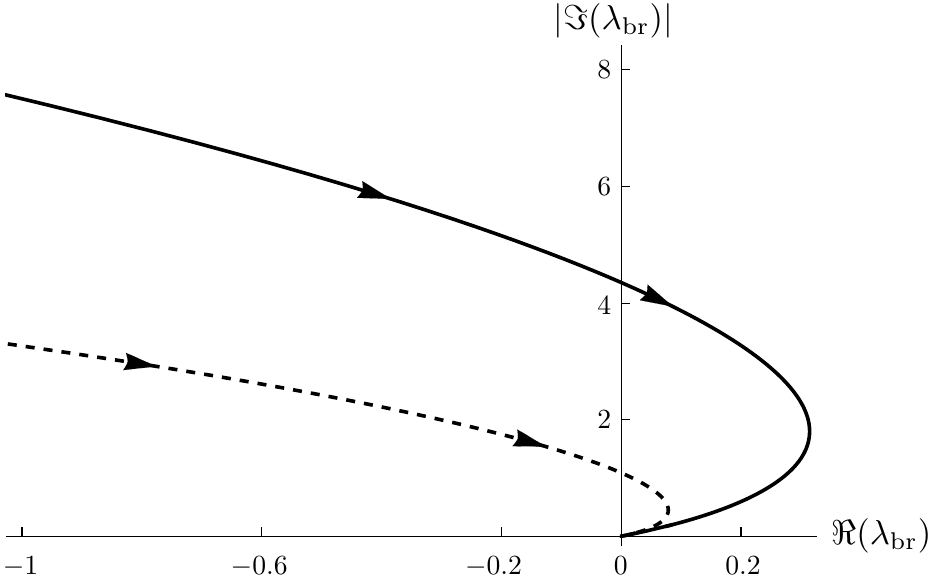}}
  \caption{Plot of the real component of the branch points versus the magnitude of the imaginary component of the branch points parametrised by $\beta>1$ for $m=0, \e=0$ and $c=1$ (dashed line) and $c=2$ (solid line). For both curves the intersections with the imaginary axis away from the origin correspond to $\beta=\beta_{\rm crit}$ and $\lim_{\beta\to\infty}|\lambda_{br}|=0$. 
Note that the figure is qualitatively similar for $0<m<1$.
}\label{FIG:beta_limit}
\end{figure}

\begin{remark}
For $0\leq m<1, \e=0$ and $\beta>\beta_{\rm crit}^m$, the absolute spectrum contains values in the right half plane. However, for a large chemotactic parameter, \ie\ $\beta\gg1$, the end points of the absolute spectrum $\lambda_{br}^\pm$ approach zero, see Figure \ref{FIG:beta_limit}.
Actually, in the limit $\beta\to\infty$, the discriminant of the characteristic polynomial of $M_m^{-}$ \eqref{EQ:disc_m} reduces to the discriminant of the characteristic polynomial of $M_1^{-}$ \eqref{EQ:discriminant_m_1}. 
That is, the branch points $\lambda_{br}^\pm$ of the absolute spectrum approach the origin from the right. 
Furthermore, the ideally weighted essential spectrum for $0\leq m<1, \e=0$ and $\beta$ large is qualitatively similar 
to the ideally weighted essential spectrum shown in the right panel of Figure~\ref{FIG:allspect_m_1} for $m=1$ and $\e=0$.
\end{remark}
\section{Small diffusion}
\label{SS:MNONZERO_EPSNONZERO}
In this section, we \edit{finish} the proof of Theorem \ref{TH:MAIN1} and show that the results obtained for $\e=0$ \edit{persist to leading order when we allow for small diffusion of the attractant $u$ in \eqref{EQ:KStw1} (\ie for $0<\e\ll1$)}.
In particular, we show that for $|\lambda|=\mathcal{O}(1)$ the weighted essential spectrum and absolute spectrum correspond, in leading order, to the spectra in the $\e=0$ case. For \edit{$|\lambda|$} large, the spectra differ significantly, however, the differences do not alter the explicit stability results since they occur in the open left half plane. 

\subsection{Set-up}
We treat the various consumption rates $0\leq m\leq1$ simultaneously. First, we eliminate the perturbation $q$, and its derivatives, from \eqref{EQ:Lop_gen}.
From the first row of \eqref{EQ:Lop_gen} we have 
\begin{align}
q  =  \, \, &\e  u^{-m} p_{zz}  +  c u^{-m}p_z - (mwu^{-1}+ \lambda u^{-m}) p. \label{EQ:q_sub_mnonzero}
\end{align}
Differentiating \eqref{EQ:q_sub_mnonzero} gives
\begin{equation}\label{EQ:qz_sub_mnonzero}
\begin{split}
q_z  = \, \, & \e  u^{-m} p^{(3)}  +\left((\e  u^{-m})_z + 
 c u^{-m}\right)p_{zz}  \\ 
& \qquad +((c u^{-m})_z-(mwu^{-1}+ \lambda u^{-m}))p_zq+ (mwu^{-1}+ \lambda u^{-m})_zp, \\ 
q_{zz}  = \, \, & \e  u^{-m} p^{(4)}+ \left(2(\e  u^{-m})_z +  c u^{-m}\right)p^{(3)}\\
\, \, & +\left((\e  u^{-m})_{zz}+2(c u^{-m})_z-(mwu^{-1}+ \lambda u^{-m})\right)p_{zz}\\
\, \, &+\left((cmu^{-m})_{zz}-2(mwu^{-1}+ \lambda u^{-m})_z\right)p_z+(mwu^{-1}+ \lambda u^{-m})_{zz} p.
\end{split}
\end{equation} 
We substitute \eqref{EQ:q_sub_mnonzero} and \eqref{EQ:qz_sub_mnonzero} into the second row of \eqref{EQ:Lop_gen} $\mL_p p+\mL_qq=\lambda q$. The resulting singular fourth ODE is
\begin{align}
\e p_{zzzz}- \cD_{m,\e}p_{zzz}- \cC_{m,\e}p_{zz}- \cB_{m,\e} p_z - \cA_{m,\e} p = 0 \label{EQ:Pop_m_nonzero_eps_nonzero}
\end{align}
where
\begin{equation*}
\begin{split}
\mathcal{A}_{m,\e}:=&(\beta+m) \left(c^2+\lambda +\lambda  m\right)\frac{ u_z^2}{u^2}-2 c \lambda  m\frac{ u_z}{u}-c (\beta+m)\frac{ u_z
   u_{zz}}{u^2}-\lambda ^2\\
   &-\lambda  (\beta+m)\frac{ u_{zz}}{u}- c (\beta -2)(\beta+m)\frac{ u_z^3}{u^3} \\
   & + \e  \bigg(c (\beta+m)\frac{ u_z u_{zz}}{u^2}-(\beta -2) (\beta+m) \frac{u_z^2 u_{zz}}{u^3}
   -(\beta+m)\frac{
   u_{zz}^2}{u^2}-\lambda  m\frac{ u_{zz}}{u}\bigg), \\
\mathcal{B}_{m,\e}:=&2 c \lambda-\left(\beta  c^2+\lambda (\beta+2 m)\right)\frac{u_z }{u}+c (\beta -m-3) (\beta+m)\frac{ u_z^2}{u^2}+c (\beta+m)\frac{ u_{zz}}{u}\\
&+\e  \left((\beta -2) (\beta+m)\frac{ u_z u_{zz}}{u^2}-c   (\beta+m)\frac{ u_{zz}}{u}\right),\\
\mathcal{C}_{m,\e}:=&-c^2+c (2 (\beta+m) +m)\frac{
   u_z}{u}+\lambda+\e \bigg(\lambda-(m+1) (\beta+m)\frac{ u_z^2}{u^2}+c m\frac{ u_z}{u}\\
   & +2 (\beta+m) \frac{u_{zz}}{u}\bigg),\\
\mathcal{D}_{m,\e}:=&-c+\e \left((\beta+2 m) \frac{u_z}{u}-c\right),
\end{split}
\end{equation*}
with $(u,w)$ the travelling wave solutions given, to leading order, by \eqref{tw_profiles_0}. 
We set $p_1:=p_z$, $p_2:=p_{zz}$ and $p_3:=p_{zzz}$ and define the operator $\mathcal{T}_\e$ by
\begin{align*}
&\mathcal{T}_{m,\e}(\lambda)\begin{pmatrix}
p\\p_1\\p_2\\p_3
\end{pmatrix}:=\begin{pmatrix}
p\\p_1\\p_2\\p_3
\end{pmatrix}'- M_{m,\e}(z,\lambda)\begin{pmatrix}
p\\p_1\\p_2\\p_3
\end{pmatrix}=0,
\end{align*}
where
\begin{align*}
M_{m,\e}(z,\lambda):=\begin{pmatrix}
0&1&0&0\\ 0&0&1&0\\0&0&0&1\\
\mathcal{A}_{m,\e}/\e&\mathcal{B}_{m,\e}/\e&\mathcal{C}_{m,\e}/\e&\mathcal{D}_{m,\e}/\e
\end{pmatrix}.\numberthis\label{EQ:Teps}
\end{align*}
All terms in $\mathcal{T}_{m,\e}$ can be expressed in terms of either $u_z/u$ or $w/u$,
since $u_{zz}=\left(cu_z-w\right)/\e$ and $w_z=-cw+\beta\left(\frac{wu_z}{u}\right)$ \eqref{EQ:KStw2}. 
Using \eqref{EQ:limits}, the limits of $\mathcal{A}_{m,\e}$, $\mathcal{B}_{m,\e}$, $\mathcal{C}_{m,\e}$ and $\mathcal{D}_{m,\e}$ as $z\rightarrow\pm\infty$ are
\begin{align*}
\mathcal{A}_{m,\e}^+:=-\lambda ^2, \qquad \mathcal{B}_{m,\e}^+:=2 c \lambda,\qquad\mathcal{C}_{m,\e}^+:=-c^2+\lambda(1 +\e),\qquad  \mathcal{D}_{m,\e}^+:=-c(1+\e),
\end{align*}
and
\begin{equation}
\begin{split}
 \cA_{m,\e} ^- = \, \, & \frac{c^4 m (\beta +m)}{(\beta +m-1)^3}-\lambda
   ^2-\frac{c^2 \lambda  m (\beta +m-2)}{(\beta +m-1)^2} - \e \left(\frac{ \lambda c^2 m }{(\beta
   +m-1)^2}-\frac{c^4 m (\beta +m)}{(\beta
   +m-1)^4} \right)\\ 
 \cB_{m,\e}^- =   \, \, & \frac{c \lambda(\beta
   -2)   }{(\beta
   +m-1)}-\frac{c^3 (\beta +m (\beta +m+2))}{(\beta
   +m-1)^2}+ \e\frac{c^3 (m+1) (\beta +m)}{(\beta +m-1)^3}  , \\ 
   \cC_{m,\e}^- =   \, \, &  \left(\frac{c^2 (\beta +2 m+1)}{\beta
   +m-1}+\lambda\right) +\e  \left( \frac{\beta  c^2}{(\beta +m-1)^2}+\lambda \right) , \\ 
   \cD_{m,\e}^- = \,\, &-c+ \e \frac{c (m+1) }{\beta +m-1} \,.
\end{split}
\end{equation}
We define the asymptotic matrices $M_{m,\e}^{\pm}(\lambda):=\displaystyle\lim_{z\rightarrow\pm\infty}M_\e^m(z,\lambda)$. That is,
\begin{align}
\label{EQ:M2}
 M_{m,\e}^{\pm}(\lambda)=\begin{pmatrix}
0&1&0&0\\0&0&1&0\\ 0&0&0&1\\
\mathcal{A}_{m,\e}^{\pm}/\e&\mathcal{B}_{m,\e}^{\pm}/\e&\mathcal{C}_{m,\e}^{\pm}/\e&\mathcal{D}_{m,\e}^{\pm}/\e
\end{pmatrix}.
\end{align}

\subsection{Proof of Theorem~\ref{TH:MAIN1} for $0 \leq m \leq 1$ and $0<\e  \ll1$}
\label{SS:PROOF2}
\edit{The matrices $M_{m,\e}^{\pm}$ have four spatial eigenvalues, while $M_{m}^{\pm}$ have only three. }
We show that the fourth spatial eigenvalue is far into the left half plane for both asymptotic matrices $M_{m,\e}^{\pm}$ (and for $|\lambda|=\mathcal{O}(1)$), while the other three spatial eigenvalues are, to leading order, given by the spatial eigenvalues of 
$M_{m}^{\pm}$. 

The characteristic polynomial of $M_{m,\e}^{+}$ is 
\begin{align}
\e\left(\mu^4+c\mu^3-\lambda\mu^2\right) +(\mu^2+c\mu-\lambda)(c\mu-\lambda)=0,\label{EQ:eps_nonzero_charpolyplus}
 \end{align}
which is regular in $\lambda$, but \edit{singularly perturbed} in $\mu$. 
In the limit $\e \to0$, we recover the characteristic polynomial of $M_m^{+}$.
The dispersion relations of $M_{m,\e}^{+}+\nu_+ I$ are
\begin{equation}\label{EQ:dispersion_plus_epsnonzero}
\lambda=-k^2-\nu_+   (c-\nu_+  )+i (c k-2 k \nu_+  )\,,\quad
\lambda=-\e k^2-\nu_+  (c-\nu_+  \e )+i (c k-2 \e k \nu_+  ) .
\end{equation}
For $\nu_+ \in (0,c)$, \eqref{EQ:dispersion_plus_epsnonzero} is fully contained in the open left half plane and the ideal weight is still $\nu_+^*=c/2$.
Observe that, unlike the $\e=0$ case, both dispersion relations of $M_\e^{m,+}+\nu_+ I$ are parabolas in $k$ and consequently
they no longer approach a vertical line in the limit $|k| \to \infty$.
 
The spatial eigenvalues of \eqref{EQ:eps_nonzero_charpolyplus} are
  \begin{align*}
  &\mu_1^+=\frac{-c+\sqrt{c^2+4\e \lambda }}{2\e}=\frac{\lambda }{c}-\frac{\lambda ^2 \e }{c^3}+\mathcal{O}(\e^2),\\
  &\mu_{2,3}^+=\frac{-c\pm\sqrt{c^2+4 \lambda }}{2}, \\
  &\mu_4^+=\frac{-c-\sqrt{c^2+4\e \lambda }}{2\e}=-\frac{c}{\e }-\frac{\lambda }{c}+\frac{\lambda ^2 \e }{c^3}+\mathcal{O}(\e^2),
  \end{align*}
where the asymptotic expansions only hold for $|\lambda|=\mathcal{O}(1)$.
 The spatial eigenvalues $\mu_{1,2,3}^+$ are, to leading order, the same as those in the $\e=0$ case  \eqref{EQ:spatial_eigenvalues_plus}. The singular spatial eigenvalue $\mu_4^+$ approaches $-\infty$ as $\e \to 0$ (for $|\lambda|=\mathcal{O}(1)$).
 
  \begin{figure}[t]
\centering
\scalebox{1}{\subfloat{\label{fig:allspect_stablea}\includegraphics{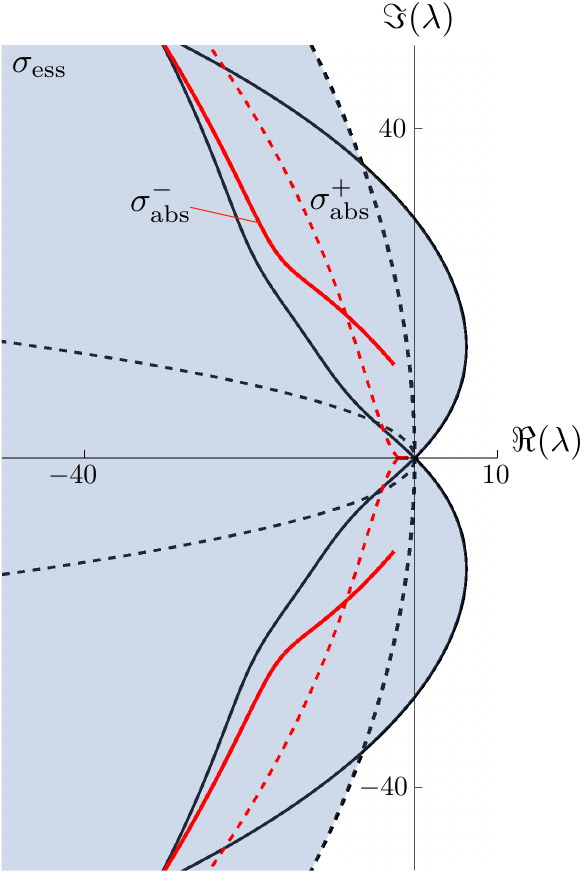}}}\hspace{.5cm}
\scalebox{1}{\subfloat{\label{fig:allspect_stableb}\includegraphics{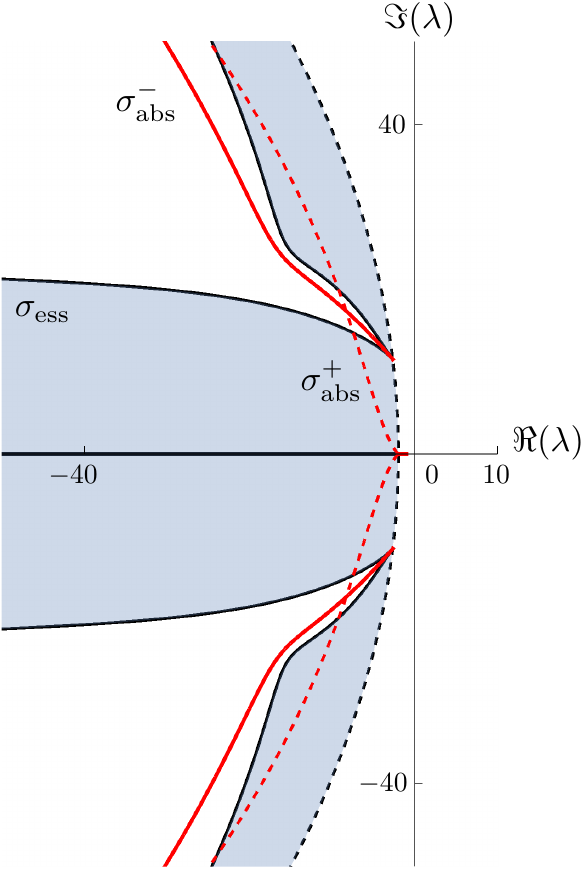}}}
\caption{The essential and absolute spectrum in the unweighted space (left panel) and in the ideally weighted space (right panel) for 
$\beta=1.3<\beta_{\rm crit}$ \eqref{EQ:10th_order_poly}, $c=2$, $\e=0.02$ and $m=0$, where the ideal weight is $\nu_-^*\approx -2.447$ and $\nu_+^*=c/2=1$.
The dispersion relations of $M_{m,\e}^{+}+\nu_+I$ \eqref{EQ:dispersion_plus_epsnonzero} are shown as black dashed lines, while those of $M_{m,\e}^{-}+\nu_-I$ 
are shown as black solid lines, $\sigma_{\rm abs}^+$ is shown as red dashed lines and $\sigma_{\rm abs}^-$ as red solid lines. The shaded regions are the interior of the (weighted) essential spectrum. 
Observe that the (weighted) essential spectra and absolute spectra agree, to leading order, for $|\lambda|=\mathcal{O}(1)$, but not for $|\lambda|$ large, to the spectra for the same parameter set but with $\e=0$, see Figure \ref{FIG:allspect_stable}. Also note that the ideal weights are similar.
}
\label{FIG:allspect_stable_eps_nonzero}
\end{figure}

The characteristic polynomial of $M_{m,\e}^{-}$ is
\begin{equation}\label{EQ:eps_nonzero_charpolyminus}
\begin{split}
&\mu ^2 \left(-\frac{c^2 (\beta +2 m+1)}{\beta +m-1}-\lambda \right)+\frac{\mu  \left(c^3 (\beta +m (\beta +m+2))-(\beta -2) c \lambda  (\beta
   +m-1)\right)}{(\beta +m-1)^2}\\
   &+\frac{c^2 \lambda  m
   (\beta +m-2)}{(\beta +m-1)^2}-\frac{c^4 m (\beta +m)}{(\beta +m-1)^3}+c \mu ^3+\lambda ^2+\e  \bigg(\frac{c^3 \mu  (m+1) (\beta +m)}{(\beta +m-1)^3}\\
   &+\mu ^2
   \left(-\frac{\beta  c^2}{(\beta +m-1)^2}-\lambda \right)+\frac{c^2 m \left(\lambda  (\beta +m-1)^2-c^2 (\beta
   +m)\right)}{(\beta +m-1)^4}\\
   &-\frac{c \mu ^3 (m+1)}{\beta +m-1}+\mu ^4\bigg).
\end{split}
\end{equation}
which is still regular in $\lambda$, but singularly perturbed in $\mu$.
In the limit $\e \to0$, we recover the characteristic polynomial of $M_m^{-}$ 
\begin{align*}
&c \mu ^3-\mu ^2 \left(\frac{c^2 (\beta+2 m+1)}{\beta+m-1}+\lambda \right)+\mu \left(\frac{ c^3 (\beta +m (r+2))}{(\beta+m-1)^2}-\frac{(\beta -2) c \lambda }{(\beta+m-1)}\right)\\
   &-\frac{c^4 m (\beta+m)}{(\beta+m-1)^3}+\frac{c^2 \lambda  m
   (\beta+m-2)}{(\beta+m-1)^2}+\lambda ^2=0,
\end{align*}
and three of the spatial eigenvalues of $M_{m,\e}^{-}$ are, to leading order, thus given by spatial eigenvalues of $M_m^{-}$ for $|\lambda|=\mathcal{O}(1)$.
We use the expansion $\mu=\eta_{-1}/\e+\eta_0+\mathcal{O}(\e)$ to determine the leading order contribution of the singular spatial eigenvalue of $M_{m,\e}^{-}$. Substituting this expansion into \eqref{EQ:eps_nonzero_charpolyminus} gives, to leading order,
$
\eta_{-1}^3(\eta_{-1}+c)=0.
$
So, the singular spatial eigenvalue of $M_{m,\e}^{-}$ is $\mu_4^- = -c/\e + \mathcal{O}(1)$ (for $|\lambda|=\mathcal{O}(1)$).
In particular, both singular spatial eigenvalues are to leading order the same and approach $-\infty$ as $\e \to 0$.

For $|\lambda|=\mathcal{O}(1)$, the (weighted) dispersion relations of $M_{m,\e}^{\pm}$ are $\mathcal{O}(\e)$ perturbations of those from $M_m^{\pm}$,
since $\mu_{1,2,3}^\pm$ are, to leading order, the same as those in the $\e=0$ case, and since the singular spatial eigenvalues $\mu_4^\pm$ have asymptotically large negative real parts (for $|\lambda|=\mathcal{O}(1)$).
Moreover, the characteristic polynomials of $M_{m,\e}^{\pm}$ are regularly perturbed in $\lambda$. 
Consequently, the Morse indices $i_{\pm}$ and the interior of the essential spectrum are unaffected 
by the singular spatial eigenvalues $\mu_4^\pm$. 
Similarly, since $\mu_4^\pm$ also does not affect the ranking of $\mu_{1,2,3}^\pm$, the absolute spectrum is, to leading order, the same as for the $\e=0$ case. In particular, the branch points $\lambda_{br}^\pm$ are, to leading order, the same as those for the $\e=0$ case and there is some parameter $\beta_{\rm crit}^m(\e)$, given to leading order by $\beta_{\rm crit}^m$, such that the branch points, and therefore the absolute spectrum, are contained in the open left half plane for $1-m<\beta<\beta_{\rm crit}^m(\e)$. 

The above asymptotic analysis is only valid for $|\lambda|=\mathcal{O}(1)$, since the singular spatial eigenvalues $\mu_4^\pm$ become $\mathcal{O}(1)$ for $|\lambda|$ large.  
However, it can be shown using asymptotic analysis that, to leading order, there are no additional intersections between the dispersion relations of $M_{m,\e}^{\pm}+\nu_\pm I$ and the imaginary axis for $|\lambda|$ large as long as $\nu_->-\frac{c(\beta+m)}{\beta+m-1}$. This condition arises from the asymptotic limits of the weighted dispersion relations $M_{m}^{-}+\nu_- I$ (see \eqref{EQ:limit22} for the analogous condition for $m=0$).
We omit the technical details of this asymptotic analysis. As the dispersion relations do not intersect the imaginary axis for large $|\lambda|$, the essential spectrum, and therefore the absolute spectrum, does not enter into the right half plane, except in the region $|\lambda|=\mathcal{O}(1)$. 
See Figure~\ref{FIG:allspect_stable_eps_nonzero} for an example of the spectral picture in the case $\e \neq 0$.
This concludes the complete proof of Theorem \ref{TH:MAIN1}.

\section{Outlook}
\label{SS:POINT}

In this manuscript, we located the weighted essential spectrum and absolute spectrum associated with travelling wave solutions to the Keller-Segel model \eqref{EQ:KStw1} for 
$0\leq m\leq 1$, $\beta>1-m$ and $0\leq \e\ll1$. 
By locating the branch points, that form the leading edge of the absolute spectrum, 
we proved that the absolute spectrum and ideally weighted essential spectrum are contained in the open left half plane for $1-m<\beta<\beta_{\rm crit}^m(\e)$
and we derived leading order expressions determining $\beta_{\rm crit}^m(\e)$.
We also developed a procedure for locating the range of weighted spaces for which the weighted essential spectrum is in the open left half plane. 
For $\beta>\beta_{\rm crit}^m(\e)$, all travelling wave solutions have absolute spectrum in the right half plane and the 
travelling wave solutions are thus absolutely unstable.
These results provide a complete picture of the absolute spectrum and weighted essential spectrum associated with all possible travelling wave solutions to the Keller-Segel model \eqref{EQ:KStw1} and they expand on the previous results for the essential spectrum known in the literature \cite{nagai1991traveling, wang2013mathematics}. 
Furthermore, it is now clear how the absolute spectrum and weighted essential spectrum deform between the limit cases $m=0$ and $m=1$. 
Moreover, we showed that the transition to the absolutely unstable parameter regime is characterised by the absolute spectrum crossing into the right half plane away from the real axis (similar to the example in \cite{rademacher2007computing}). 

In order to complete the full spectral picture for travelling wave solutions to \eqref{EQ:KStw1} the point spectrum must also be located. 
This is usually far more involved as it requires information about the linearised system on the whole spatial domain, rather than just its asymptotic behaviour.
If there exists \edit{point spectrum} with positive real part, then the travelling wave solutions are spectrally unstable, regardless of the stability properties of the absolute spectrum and weighted essential spectrum. 
Note that point spectrum is unaffected by weighting the space \cite{kapitula2013spectral}.
An early proof offered by \cite{rosen1975stability} shows that there are no positive eigenvalues for $0\leq\e\ll1$ under the assumption that eigenvalues are real-valued.
However, it is unclear that this assumption holds, since the linearised operator $\mL$ \eqref{EQ:Lop_gen} is not self-adjoint. 

An analytic tool for locating the point spectrum is the Evans function \cite{evans1,evans2,evans3,evans4}. 
Unfortunately, the Evans function is generically hard to explicitly compute for systems of partial differential equations and one has to rely on numerics.
This is also the case here, especially since the explicit travelling wave profiles for \eqref{EQ:KStw1} with $\e >0$ are not known.
In \cite{harley2015numerical}, the Evans function associated with travelling wave solutions to \eqref{EQ:KStw1} with $m=0$ and $\e=0$ was calculated numerically using a Riccati transformation. It was shown that there is a second order temporal eigenvalue at the origin and that there are no other eigenvalues in the right half plane with $|\lambda|<10^7$.
Due to the translation invariance, $\lambda=0$ persists as an eigenvalue (with order at least one) for $0<\e\ll1$. 
However, the second eigenvalue most likely perturbs for $0<\e\ll1$ determining the fate of the spectral stability of the travelling wave solution (assuming the weighted essential spectrum is in the open left half plane).   
In ongoing research, we are addressing the issue of the point spectrum by using methods similar to the ones used in \cite{harley2015numerical}.

If there is no point spectrum in the right half plane, one can conclude that the travelling wave solutions are spectrally stable in the ideally weighted space for $1-m<\beta<\beta_{\rm crit}^m(\e)$, \ie\ transiently unstable. 
Ideally, one would like to use this spectral stability result to conclude nonlinear (in)stability of the travelling wave solutions. 
For a sectorial semilinear operator with a spectral gap (\ie\ the spectrum is contained in the open left half plane except for the translation invariance eigenvalue at the origin), spectral stability implies nonlinear stability of the associated travelling wave solution \cite{henry1981geometric,sandstede2000spectral}.
However, while the operator $\mL$ \eqref{EQ:Lop_gen} appears to be sectorial for $0<\e\ll1$, see, for instance, Figure~\ref{FIG:allspect_stable_eps_nonzero}, it is 
quasilinear rather than semilinear.
In \cite{meyries2014quasi}, it was shown that for a large class of quasilinear parabolic reaction-diffusion systems one can still deduce nonlinear stability results from the spectral stability results, as long as the linearised operator fulfills certain conditions.
Unfortunately, the Keller-Segel model studied in this manuscript does not fall into the class of quasilinear parabolic reaction-diffusion systems considered in \cite{meyries2014quasi}.
For the Keller-Segel model \eqref{EQ:KSgen} with nonlinear diffusion and with logarithmic chemosensitivty (\ie\ $\Phi(u)=\log(u)$), linear consumption (\ie\ $m=1$) and nonzero growth (\ie\ $\kappa>0$),
the general theory for semilinear operators was extended in \cite{meyries2011local} to prove nonlinear instability results in certain cases of the model.
Another approach using a Hopf-Cole transformation, in conjunction with weighted energy estimates, was used in \cite{li2014stability,li2012steadily} to deduce nonlinear stability results for the Keller-Segel model \eqref{EQ:KSgen} with logarithmic chemosensitivty, linear consumption and zero growth. 
Alternatively, in order to apply the general theory
for semilinear systems, \cite{henry1981geometric} proposes to transform a quasilinear system to a semilinear system.
Observe that this approach is akin to the method used in \S\ref{S:SUBLIN}. 
It is a challenge to see if any of these methods can be used to obtain nonlinear stability results for the travelling wave solutions of \eqref{EQ:KStw1} studied in this manuscript.

The dynamical implications of absolute spectrum in the right half plane for travelling wave solutions of the Keller-Segel model \eqref{EQ:KStw1} are not known. 
In typical examples, such as the F-KPP equation, the transition to an absolutely unstable regime is associated with the so-called linear spreading speed, \ie\ the speed `generic' initial conditions will eventually travel at. Note that in the F-KPP equation this is known as the minimal wave speed. In other words, the linear spreading speed is the speed of a travelling wave solution `selected' by the model. However, in the Keller-Segel model \eqref{EQ:KStw1} the transition to the absolutely unstable regime is, to leading order, independent of the wave speed and it thus does not seem to have an influence on the asymptotic speed of a generic initial condition (that evolves to a travelling wave solution).
Rather, the initial condition of the bacteria population $w$ determines the wave speed \cite{nagai1991traveling}.
Note that in the case of a Keller-Segel model \eqref{EQ:KSgen} with a growth term, 
the absolute spectrum does appear to have an influence on the
wave speed selection \cite{bose2013invasion,meyries2011local}.
Moreover, as the transition of the absolute spectrum into the right half plane is complex valued, one expects oscillatory instabilities to manifest themselves. 
These type of bifurcations have, for instance, been studied in \cite{sandstede1999essential,sandstede2001structure}. 
Future work will address, both analytically and numerically, what the absolute instabilities imply dynamically and the connection, if any, with the wave speed.

\section*{Acknowledgments}
PD and
PvH acknowledge support under the Australian Research Council's Discovery Early Career Researcher Award funding scheme DE140100741. \edit{RM would like to thank M. Holzer and J. Rademacher for very informative discussions regarding the material covered in this paper as well as for pointing out some illuminating references}.

\bibliographystyle{siamplain}
\bibliography{ABSOLUTE_INSTABILITIES_KS_AUG_2016}
\end{document}